\documentclass{article}
\usepackage{amsmath, amsthm, amssymb, enumerate}

\newcounter{test}
\newcounter{casest}[paragraph]

\newtheorem{theorem}{Theorem}[section]
\newtheorem{corollary}[theorem]{Corollary}
\newtheorem{proposition}[theorem]{Proposition}
\newtheorem{lemma}[theorem]{Lemma}

\theoremstyle{definition}
\newtheorem{definition}[theorem]{Definition}
\newtheorem{remark}[theorem]{Remark}
\newtheorem{example}[theorem]{Example}
\newtheorem{cs}[casest]{Case}
\newcommand{\g}{\mathfrak{g}}
\newcommand{\h}{\mathfrak{h}}

\numberwithin{equation}{section}

\begin{document}
\title{Classification of Lie bialgebras over current algebras}
\author{F.~Montaner, A.~Stolin, and E.~Zelmanov}
\date{}

\maketitle



\thispagestyle{empty}

\begin{abstract}
In this paper we give a classification of Lie bialgebra structures
on Lie algebras of type $\g[[x]]$ and $\g[x]$, where $\g$ is a simple
complex finite dimensional Lie algebra.

\smallskip

{\bf Mathematics Subject Classifications (2010).} 17B37, 17B62.

{\bf Key words:} Lie bialgebras, quantum groups, orders in loop algebras
\end{abstract}

\section{Introduction}

In what follows $F$ is an algebraically closed field of characteristic zero. By a quantum group we mean a Hopf algebra $A$ over the series $F[[h]]$, where $h$ is a formal parameter such that:

1. ${A}/{hA}\cong U(L)$ as a Hopf algebra. Here $U(L)$ is the universal
enveloping algebra of some Lie algebra $L$ viewed as a Hopf algebra with
comultiplication $\Delta_{0}=a\otimes 1+1\otimes a,\ a\in L$.

2. $A\cong V[[h]]$ as a topological $F[[h]]$-module for some vector space
$V$ over $F$.

The important example of quantum groups are the quantum universal enveloping
algebra $U_{h}(\g)$, the quantum affine Kac-Moody algebra
$U_{h}(\widehat{\g})$ and the Yangian $Y(\g)$, here $\g$ in a
finite-dimensional simple complex Lie algebra.

If $A$ is a quantum group and ${A/hA}\cong U(L)$ then $L$ is equipped with a Lie bialgebra structure
$\delta: L\longrightarrow L\otimes L$, $\delta(b)=h^{-1}(\Delta(a)-\Delta^{op}(a))$ mod$\ (h)$, $b\in L$, here $a$ is a preimage of $b$ in $A$.
The Lie bialgebra $(L,\delta)$ is called the classical limit of the quantum group
$A$.

P. Etingof and D. Kazhdan proved that an arbitrary Lie bialgebra is a classical limit of some quantum group (see \cite{EK, EK1}).
Moreover, they proved that for any Lie algebra $L$ there exists a one-to-one correspondence between Lie
bialgebra structures on $L[[h]]$ and those quantum groups which have $L$ with some $\delta$ as the classical limit.

The aim of this paper is classification of Lie bialgebras on the current
algebras $\g[[x]]$ and $\g[x]$, which is an important step towards
classification of quantum groups such that ${A/hA}\cong U(\g[[x]])$ or
${A/hA}\cong U(\g[x])$.

We proceed in the following way.

I. Any Lie bialgebra structure $\delta$ on a Lie algebra $L$ gives rise to an embedding $L\subset D(L,\delta)$ into the Drinfeld double algebra $D(L,\delta)$
equipped with a nondegenerate skew-symmetric bilinear form (see \cite{DR}).

We prove in Theorem \ref{thm:1} that there are four types of possible
Drinfeld doubles of bialgebras on $\g[[x]]$. They correspond to the
classical $r$-matrices
$$
r_{1}(x,y)=0,\ r_{2}(x,y)=\frac{\Omega}{x-y},\
r_{3}(x,y)=\frac{x\Omega}{x-y}+r_{DJ},\ r_{4}(x,y)=\frac{xy\Omega}{x-y}
$$
where $\Omega$ is a Casimir (invariant) element of the $\g$-module
$\g\otimes \g$ and $r_{DJ}$ is the Drinfeld-Jimbo classical $r$-matrix which
will be defined below.

Throughout this paper we use the following root space decomposition (or
Cartan decomposition) $\g=\h\oplus (\oplus_\alpha\g_\alpha)$. Here $\h$ is a
Cartan subalgebra and $\{\alpha\}$ is the set of roots with respect to $\h$.
We also denote by $\alpha_i$, $i=1,2,\ldots,\rm{rank}(\g)$, the set of
simple roots, and $\alpha_0=-\alpha_{\max}$. We define the positive integers
$k_i$ using the following relation:
$\sum_{i=0}^{\rm{rank}(\g)}k_i\alpha_i=0$, where $k_0=1$.

The Drinfeld-Jimbo classical $r$-matrix is defined by the following formula:
$$
r_{DJ}=\sum_{\alpha>0}e_\alpha\otimes e_{-\alpha}
+\frac{1}{2}\sum_{i=1}^{\rm{rank}(\g)}h_{\alpha_i}\otimes h'_{\alpha_i}
$$
Here $e_{\alpha}\in\g_\alpha$ are chosen such that the Killing form
$K(e_\alpha ,e_{-\alpha})=1, h_{\alpha_i}=[e_{\alpha_i},e_{-\alpha_i}]$ and
$\{h'_{\alpha_i}\}$ is the basis of $\h$ dual to $\{h_{\alpha_i}\}$ with respect to the Killing form.

We also consider the multivariable case $\g[[x_1,x_2,\ldots,x_n]]$. We prove
in Theorem \ref{thm:4} that there exists only the trivial double over
$\g[[x_1,x_2,\ldots,x_n]]$ for $n\geq 2$.

II. Using the classification of the Drinfeld doubles on $\g[[x]]$, in
Theorems \ref{theorem-ma}, \ref{thm-B}, \ref{theorem-C} we obtain a
classification of the doubles on $\g[x]$. The classification table contains
seven separate cases and a family of algebras, which can be parametrized by
$\mathbb{CP}^1/\mathbb{Z}_2$.

III. If Lie bialgebras $(L,\delta_{1})$ and $(L,\delta_{2})$ lead to the
same Drinfeld double, then $\delta_{1}-\delta_{2}$ is called the classical
twist. Using the theory of the maximal orders we classify classical twists
in all the cases for $\g[x]$ (see Theorems \ref{theorem-D}, \ref{theorem-E},
\ref{theorem-F}, \ref{theorem-G}).

\thanks{The authors are grateful to the Korean Institute for Advanced Study
for the hospitality during the preparation of the manuscript. The first
author was supported by the Spanish Ministerio de Educaci\'on y Ciencia and
FEDER (MTIM 2007-67884-C04-02) and Diputaci\'on General de Arag\'on (Grupo
de Investigaci\'on de \'Algebra). The second author was supported by the
Crafoord Foundation and the Swedish Research Council. The third author was
partially supported by the NSF.}

\section{Classification over series}

In what follows $L = \g[[x]]$ is the algebra of infinite series over a
finite dimensional simple Lie algebra $\g$. Let $M$ be the maximal ideal of
$F[[x]]$ which consists of series with zero constant term. A functional
$f:L\rightarrow F$ (resp. $F[[x]] \rightarrow F$) is called a distribution
if there exists $m \ge 1$, such that $f(\g \otimes M^m) =0$ (resp. $f(M^m) =
(0)$). Let $L^*$ denote the space of all distributions on $L$. Let $L
\otimes L$ be the (topological) tensor square of $L$. A comultiplication
$\delta: L \rightarrow L \otimes L $ defines an algebra structure on $L^*$
via \[ \delta(a) = \sum a_{(1)} \otimes a_{(2)}; f,g \in L^*, (fg)(a) = \sum
f (a_{(1)})g(a_{(2)}).\]

Following V. Drinfeld (see \cite{DR}) we call $(L, \delta)$ a Lie bialgebra if
\begin{itemize}
\itemsep-3pt
\item
the induced algebra structure on $L^*$ is a Lie algebra,
\item
$\delta$ is a derivation of $L$ into the $L$-module $L \otimes L$.
\end{itemize}
The vector space $D = L \oplus L^*$ is equipped with a natural nondegenerate symmetric bilinear form
$$
(L|L)= (L^*|L^*) = (0),\quad
(a|f) = f(a)
$$
for $a \in L, f \in L^*$.  Moreover, there exists a unique Lie algebra structure on $D$, which extends the multiplications on $L$ and $L^*$ respectively and makes the form defined above invariant. We call $D$ the Drinfeld double of the bialgebra
$(L,\delta)$ and denote it $D(L,\delta)$.

\begin{lemma}
\label{lemma:1} Let $L=\g[[x]]$. As a $\g$-module $D(L,\delta)$ is a direct
sum of regular (adjoint) $\g$-modules. 
\end{lemma}

\begin{proof}
Let $U(\g)$ be the universal enveloping algebra of the Lie algebra $\g$. For a $\g$-module $V$ let $\mathrm{Ann}(V)$ denote the ideal $\mathrm{Ann}(V) = \{ a \in U:\ aV =(0)\}$ We will start with the following observations: if $V,W$ are finite dimensional irreducible $\g$-modules and $\mathrm{Ann}(W)V=(0)$ then $V \cong W$. Indeed, let $w_1, \dots, w_n$ be a basis of $W$, $\mathrm{Ann}(w_i) = \{a \in U:\ a w_i=0\}, 1 \le i \le n$, are the corresponding maximal left ideals of $U(\g)$, \[\mathrm{Ann}(W)= \mathrm{Ann}(w_1) \cap \dots \cap \mathrm{Ann}(w_n). \]

If $v$ is a nonzero element of $V$ then by our assumption $\mathrm{Ann}(W)
\subset \mathrm{Ann}(V)$. Hence $V \cong U(\g) / \mathrm{Ann}(v)$ is a
homomorphic image of the module $U(\g)/\mathrm{Ann}(W)$, which is a
submodule of $U(\g)/\mathrm{Ann}(w_1) \oplus \dots \oplus U(\g)/
\mathrm{Ann}(w_n) \cong W \oplus \dots \oplus W$. This implies that $V \cong
W$.

Let $P = \mathrm{Ann}(\g)$ be the annihilator of the regular $\g$-module.
Clearly $PL=(0)$. The algebra $U(\g)$ has an involution $\ast:U(\g)
\rightarrow U(\g)$, which sends an arbitrary element $a \in \g$ to $-a$.
Since the algebra $\g$ is equipped with a symmetric nondegenerate bilinear
form it follows that $P^*=P$.

For an arbitrary $x \in U(\g)$, arbitrary element $d_1, d_2 \in D$ we have
$(x d_1 | d_2) = (d_1 | x^* d_2)$. Hence \[ (PD|L) = (D|P^*L) = (D|PL)=(0).
\]

Hence, $PD \subset L^\perp = L$. Therefore $PPD =(0).$ Since $P$ is the kernel of the homomorphism $U(\g) \rightarrow \mathrm{End}_F(\g)$, it has finite codimension in $U(\g)$. Let
$$
U(\g)= P + \sum_{i=1}^k F a_i, a_i \in U(\g).
$$
Since the algebra $U(\g)$ is Noetherian the ideal $P$ is finitely generated as a left ideal, $P = \sum_{i=1}^l U(\g) p_i$. Now, for an arbitrary element $d \in D$ we have \begin{align*} U(\g)d & = \sum_{i=1}^k Fa_i d + Pd \\ & = \sum_{i=1}^k F a_i d + \sum_{j=1}^l U(\g) p_j d \\ & \subseteq \sum_i F a_i d + \sum_{i,j} F a_i p_j d + P^2 d \\ & = \sum_i F a_i d + \sum_{i,j} F a_i p_j d. \end{align*}

We proved that $\dim_F U(\g) d < \infty$ and therefore $D$ is a direct sum
of irreducible finite dimensional $\g$-modules. From the observation above
and from $P^2D=(0)$ it follows that all these irreducible modules are
isomorphic to the regular module. The lemma is proved. \end{proof}

\begin{corollary} There exists an associative commutative $F$-algebra $A$ with 1 such that $D \cong \g \otimes_F A$.
\end{corollary}

\begin{proof}
The corollary immediately follows from the lemma above and from Proposition
2.2 in \cite{BZ}.
\end{proof}

\begin{lemma}
\label{lemma:2}
Let $A$ be an associate commutative $F$-algebra with 1 and let $\g$ be a simple finite dimensional Lie algebra over $F$. Suppose that $( \ | \ )$ is a symmetric invariant nondegenerate bilinear form on $\g \otimes_F A$. Then there exists a linear functional $t: A \rightarrow F$ such that for arbitrary elements $x,y \in \g$ and $a,b \in A$ we have
$$
(x \otimes a | y \otimes b) = K(x,y) \cdot t(ab),
$$
where $K(x,y)$ is the Killing form on $\g$.
\end{lemma}

\begin{proof}
Consider the root decomposition $\g = \h\oplus(\oplus_{\alpha} \g_\alpha)$.
Fix $\alpha \in \Delta$ and choose elements $e_\alpha \in \g_\alpha,
e_{-\alpha} \in \g_{-\alpha}$ such that $K(e_\alpha,e_{-\alpha}) =1$. Define
the functional $t: A \rightarrow F$ via $t(a)=(e_\alpha\otimes 1|\
e_{-\alpha}\otimes a)$.

Choose a basis $e_1, \dots, e_n$ in $\g$ and the dual basis $e^1, \dots, e^n$ with respect to the Killing form. Suppose that the Casimir operator $C = \sum_{i=1}^n a d(e_i)a d(e^i)$ acts on $\g$ as multiplication by $\gamma \in F, \gamma \ne 0$. Choose $a,b \in A$ and let
$C(a): \g\otimes A\to \g\otimes A$ be defined by
$$
C(a) = \sum_{i=1}^n ad (e^i\otimes 1) ad (e_i \otimes a).
$$
Then $e_\alpha \otimes ab = \frac{1}{\gamma} C(a) (e_\alpha\otimes b)$.
Since the form $(\ |\ )$ is invariant, we have

\begin{align*}
( e_\alpha \otimes a | e_{-\alpha} \otimes b)
& = \frac{1}{\gamma}(C(a) (e_\alpha\otimes 1) | e_{-\alpha} \otimes b ) \\
&=  \frac{1}{\gamma} ( e_\alpha\otimes 1 | C(a) (e_{-\alpha} \otimes b) )\\
&= (e_\alpha\otimes 1 | e_{-\alpha} \otimes ab) \\
& = t (ab)
\end{align*}
Now let $\beta$ be an arbitrary root. Choose $e_\beta \in \g_\beta, e_{-\beta} \in \g_{-\beta}$. We have
$$
\left([ e_\alpha, e_\beta] \otimes a\ |\ [ e_{-\alpha}, e_{-\beta}] \otimes b\right)
= (e_\alpha\otimes 1 |\  [ e_\beta, [e_{-\alpha}, e_{-\beta}]] \otimes ab).
$$
Let $[e_\beta, [e_{-\alpha}, e_{-\beta}]] = \xi e_{-\alpha}, \xi \in F$.  Then
\[
(e_\alpha\otimes 1 |\ [e_\beta, [e_{-\alpha},e_{-\beta}]] \otimes ab) =
\xi(e_\alpha \otimes 1| e_{-\alpha}\otimes ab) = \xi t(ab).
\]
On the other hand,
\[
K ( [ e_\alpha, e_{\beta}],[e_{-\alpha}, e_{-\beta}]) = K (
e_\alpha,[e_\beta, [e_{-\alpha}, e_{-\beta}]]) = \xi.
\]
Hence,
$$
([e_\alpha, e_\beta] \otimes a |\  [e_{-\alpha}, e_{-\beta}] \otimes b)
= K ( [e_\alpha, e_\beta] , [e_{-\alpha}, e_{-\beta}]) t(ab).
$$
This implies the lemma.
\end{proof}

Let $A$ be an associative commutative unital $F$-algebra with a functional
$t: A \rightarrow F$.

\begin{definition}
We call $(A,t:A\rightarrow F)$ a trace extension of $F[[x]]$ if
\begin{enumerate}[(i)]
\item the bilinear form $(a|b) = t(ab)$ is nondegenerate,
\item $F[[x]]^\perp = \{a \in A:\ (F[[x]]|a) = (0) \} = F[[x]]$,
\item for an arbitrary distribution $f: F[[x]] \rightarrow F$ there exists an element $a \in A$ such that $f(b) = (b|a)$ for an arbitrary $b \in F[[x]]$.
\end{enumerate}

\end{definition}

\begin{remark}
\label{rem:1} If $F[[x]] \subset A' \subset A''$, both $A'$ and $A''$ are
trace extensions of $F[[x]]$ and the trace of $A''$ is an extension of the
trace of $A'$, then $A' = A''$.\end{remark}

\begin{example} Let $F((x_e))$ be the algebra of Laurent series in $x_e$ with the identity $e$.  Let
$$
Ff + Fx_f + \dots + F x_f^{n-1} = F[x_f | x_f^n=0]
$$
be the algebra of truncated polynomials in $x_f$ with the identity $f; n \ge 1$. The algebra $F[[x]]$ embeds into the direct sum $\left ( \sum_{i=0}^{n-1} F x_f^i\right) \oplus F((x_e))$ via $x \rightarrow x_f +x_e.$ Choose a sequence $\alpha = (\alpha_i \in F, -\infty < i \le n-2)$ and define a trace
$$
t: \left ( \sum_{i=0}^{n-1} F x_f^i\right) \oplus F((x_e)) \rightarrow F
$$
via $t(x_e^i) =0$ for $i \ge n$, $t(x_e^{n-1}) =1$, $t(x_e^i) = \alpha_i$ for $i \le n-2$ and $t(x_f^{n-1}) = -1$, $t(x_f^i) = -\alpha_i$ for $0 \le i \le n-2$ if $n \ge 2$.

The algebra
$$
A(n,\alpha) = \left( \sum_{i=0}^{n-1}Fx_f^i \right )\oplus F((x_e))
$$
with the trace $t$ is a trace extension of $F[[x]]$. If $n=0$ then we let $A(n,\alpha) = F((x)), t(x^{-1})=1$, $t(x^i)=\alpha_i, i \le -2$.

\end{example}

\begin{example} (trivial extension) Let $A(\infty) = \sum_{i=0}^\infty F a_i + F[[x]]$ with the multiplication $a_ia_j =0, a_i x^j = a_{i-j}$ if $i \ge j$, otherwise $a_i x^j =0; t(a_0)=1, t(a_i) =0$ for $i \ge 1$.
\end{example}

\begin{definition} Two trace extensions $(A,t)$ and $(A',t')$ are isomorphic if there exists an algebra isomorphism $f:A\rightarrow A'$ which is identical on $F[[x]]$ and a nonzero scalar $\xi \in F$ such that $t'(f(a)) = \xi t(a)$ for an arbitrary $a \in A$.
\end{definition}

Let $(A,t)$ be a trace extension of $F[[x]]$ and let $\phi$ be an automorphism of the algebra $A$, such that $F[[x]]^\phi = F[[x]]$. Define a new trace $t^\phi:A\rightarrow F$ via $t^\phi(a) = t(a^\phi)$. Then $A^{(\phi)} = (A, t^\phi)$ in a trace extension of $F[[x]]$ as well.

Consider the group $G = \mathrm{Aut}_0F[[x]]$ of infinite series $x + \gamma_2 x^2 + \gamma_3 x^3 + \dots, \gamma_i \in F$, with respect to substitution. For $\phi \in G$ the mapping $x_e \rightarrow \phi(x_e), x_f \rightarrow \phi(x_f)$ induces an automorphism of the algebra $A(n,\alpha), F[[x]]^\phi = F[[x]]$. Clearly there exists a sequence $\beta = (\beta_i, i \le n-2)$ such that $A(n,\alpha)^{(\phi)} = A(n,\beta)$.

\setcounter{test}{1}

\begin{proposition}
Let $(A,t)$ be a trace extension of $F[[x]]$. Then either $(A,t) \cong A(\infty)$ or $(A,t) \cong A(n,\alpha)$ for some $n \ge 0$,
$\alpha = (\alpha_i \in F,\ -\infty < i \le n-2)$.
\end{proposition}

\begin{proof}
Consider the descending chain of ideals $A \ge xA \ge x^2 A \ge \dots$

\begin{cs} All inclusions $x^n A > x^{n+1}A, n \ge 0$ are strict. We claim that in this case $x^nA \cap F[[x]] = x^n F[[x]]$ for an arbitrary $n \ge 0$. For $n = 0$ the equality is obvious. Choose a minimal $n \ge 1$ such that $x^n F[[x]]\lneqq x^n A \cap F[[x]]$. Then $ x^{n} F[[x]] \lneqq x^n A \cap F[[x]] < x^{n-1} A \cap F[[x]] = x^{n-1}F[[x]]$. The codimension of $x^nF[[x]]$ in $x^{n-1}F[[x]]$ is 1. Hence $x^nA \cap F[[x]] = x^{n-1} A \cap F[[x]]$. This implies that $x^{n-1} \in x^n A$ and therefore $x^{n-1}A   \subset x^n A$, which contradicts our assumption. Since $F[[x]] \subset A$ is a trace extension there exists an element $a_0' \in A$ such that $t(a_0')=1$, $t(Ma_0') =0$. Then $xa_0' \in F[[x]]^\perp = F[[x]]$. By the above $xa_0' \in x F[[x]]$. Hence there exists an element $b_0 \in F[[x]]$ such that $xa_0' = xb_0$. Let $a_0 = a_0' -b_0$. Then $a_0 x =0$, $t(a_0)=1, t(a_0M) = (0)$. Suppose that we have found elements $a_0, \dots, a_n \in A$ such that
$$
x^i a_j =  \begin{cases}
a_{j-i} & \mathrm{for} \ j \ge i \\
0 & \mathrm{otherwise} \end{cases}.
$$
As above there exists an element $a_{n+1}' \in A$ such that $t(x^i a_{n+1}') = \delta_{i,n+1}$. Then $t(x^i(xa_{n+1}')) = \delta_{i,n}$, which implies that $xa_{n+1}' = a_n + b_n, b_n \in F[[x]]$. Now $x^{n+2} a_{n+1}' = x^{n+1} a_n + x^{n+1}b_n$. Hence $x^{n+1} b_n \in x^{n+2} A \cap F[[x]] = x^{n+2} F[[x]]$. Hence $b_n = x c_n, c_n \in F[[x]]$. Let $a_{n+1} = a_{n+1}' - c_n$. We have $x  a_{n+1} = x a_{n+1}' - b_n = a_n$. Thus all assumptions about the element $a_{n+1}$ are satisfied.

For arbitrary $i,j \ge 0$ we have $a_j = x^{i+1} a_{i+j+1}$, which implies $a_i a_j  = (x^{i+1}a_j)a_{i+j+1} =0$. Now $A$ contains a subalgebra $A(\infty) = \mathrm{{\rm Span}}(a_{i,j}, i \ge 0) + F[[x]]$ which is a trace extension of $F[[x]]$. In view of the Remark \ref{rem:1} $A = A(\infty)$.
\end{cs}

\begin{cs}
There exists $n \ge 0$ such that $x^n \in x^{n+1} A$. Let $n$ be minimal with this property. Let $x^n = x^{n+1}a, a \in A$. Consider the element $e = x^n a^n$. We have $e^2 = (x^{2n} a^n) a^n = x^n a^n =e$.

If $n =0$ then $e =1$.
Let $f = 1-e$ and denote $x_e = xe, x_f = xf$. The element $x_e$ is invertible in $$
eA: x_e^{-1}  =
\begin{cases}
x^{n-1}a^n & \mathrm{if} \quad n\ge1 \\
a & \mathrm{if} \quad n =0
\end{cases}.
$$
From $ax^{n+1} = x^n$ it follows that $e x^n = a^n x^{2k} = x^n$.
Hence $x_f^n =0$.
By the minimality of $n$, $ax^n \ne x^{n-1}$, which implies that $e x^{n-1} \ne x^{n-1}$ and therefore $x_f^{n-1} \ne 0$.

The mapping $eF[[x]] \rightarrow F((x_e))$ is an isomorphism, we will identify $eF[[x]]$ with $F((x_e))$. Now $A$ contains the subalgebra $\left ( \sum_{i=0}^{n-1} Fx_f^i\right) \oplus F((x_e))$. For an arbitrary $i \ge 0$ we have $t(x^i) = t(x_f^i) + t(x_e^i)= 0$. If $t(x_e^{n-1}) =0$ then $x_e^{n-1} \in F[[x]]^\perp = F[[x]]$, which contradicts $x_f^{n-1} \ne 0$. Hence for $n \ge 1$, $t(x_e^{n-1}) \ne 0$. It is clear that for $n=0$ the element $x^{-1}$ does not lie in $F[[x]]$ as well, and therefore $t(x^{-1}) \ne 0$. Scaling the trace we can assume that $t(x_e^{n-1}) = 1$, $t(x_f^{n-1}) = -1$. For $- \infty < i \le n-2$ let $\alpha_i = t(x_e^i)$,
 $\alpha= (\alpha_i, -\infty < i \le n-2)$. Then
 $$
 \left ( \sum_{i=0}^{n-1} F x_f^i \right) \oplus F((x_e)) \cong A(n,\alpha).$$
 By Remark \ref{rem:1}
 $$
 \left ( \sum_{i=0}^{n-1} F x_f^i \right) \oplus F((x_e))  = A.
$$
The proposition is proved.
\end{cs}
\end{proof}

Let $\g$ be, as above, a simple finite dimensional Lie algebra,
\hbox{$L=\g \otimes_F A(n,\alpha)$}, the nondegenerate symmetric bilinnear form on $L$ is defined via   $(x \otimes u | y \otimes y) = K(x,y) t(xy)$,
where $x,y \in \g$ and $u,v \in A(n,\alpha), K(x,y)$ is the Killing form.

Our aim is to prove the following theorem.

\begin{theorem}
\label{thm:1} (i) $D$ contains a 
subalgebra $W$ such that $(W|W)
= (0)$ and $D = \g \otimes F[[x]] +W$ is a direct sum if and only if $0\leq
n \leq 2$,

(ii) For $0\leq n \leq 2$ there exists a unique $\phi \in \mathrm{Aut}_0
F[[x]]$ such that $A(n,\alpha) \cong A(n,0)^{(\phi)}$.
\end{theorem}

This theorem immediately implies that the Drinfeld double of a Lie bialgebra $(\g \otimes F[[x]],\delta)$ is isomorphic to to $\g\otimes A(\infty)$ or to
$\g \otimes A(n,0)^{(\phi)}$, $0 \le n \le 2$.

In what follows we will often write $af$ instead of $a \otimes f$ for $a \in \g, f \in A(n,\alpha)$.

Suppose that $D$ contains a 
subalgebra $W$ such that $(W|W) = (0)$ and $D = \g \otimes F[x]+W$ is a direct sum of vector spaces. Suppose that $n \ge 3$. For an arbitrary element $a \in \g$ fix elements $p(a) = \sum_{i \ge 0}p_i(a) x^i$, $q(a) = \sum_{i \ge0} q_i (a) x^i$; $p_i(a), q_i(a) \in \g$, such that $ax_e^{n-1} + p(a), a x_e^{n-2}+q(a) \in W$.

Since $W$ is isotropic, for arbitrary elements $a,b \in \g$, we have
\begin{align*}
(ax_e^{n-1} + p(a) | b x_e^{n-1} + p(b))
&= ( a x_e^{n-1} | p(b)) + (p(a)| bx_e^{n-1})\\
&= K(a,p_0(b)) + K(p_0(a),b) = 0
\end{align*}
and similarly
\[
(a x_e^{n-1} + p(a) | b x_e^{n-2} + q(b) ) = K (a,q_0(b)) + K(p_1(a),b) = 0
\]
Define $p_0^{T}: \g\to \g$ by $K(p_0(a),b)=K(a,p_0^{T}(b))$.
Then we see that the linear transformation $p_0 = -p_0^T$ is skew-symmetric with respect to $K$ and $p_1^T + q_0 =0$.
Let $\mathrm{im}\ p_0, \mathrm{im}\ q_0$ denote the images of $p_0, q_0$ respectively. Then
$$
\ker p_0 = (\mathrm{im}\ p_0)^\perp, \quad
\ker p_1 = (\mathrm{im}\ q_0)^\perp.
$$
Consider the subspace
$$
S = \ker p_0 \cap \ker p_1 = (\mathrm{im}\ p_0 + \mathrm{im} \ q_0)^\perp.
$$
In what follows $O(k)$ will stand for an element lying in $\sum_{i \ge k} \g \otimes x_f^i + \sum_{i \ge k} \g \otimes x_e^i$.

\begin{lemma}
\label{lemma:3}
If $ax_e^k + bx_f^k + O(k+1) \in W$, $a, b \in \g$, $k >1$, then $a=0$.
\end{lemma}

\begin{proof}
If $ax_e^k+bx_f^k+O(k+1)\in W$, then $k \le n-1$. Indeed, if $k \ge n$ then $x_f^k=0$, $x_e^k=x^k$ and the element lies in $\g \otimes F[[x]]$. Now choose a counterexample with a maximal $k$. Since $a \ne 0$ there exists an element $c \in \g$ such that $[[c,a],a] \ne 0$. The subalgebra $W$ contains an element $cx_e^{-1} + O(0)$, $O(0) \in \g \otimes F[[x]]$. Now,
\begin{align*}
[[c x_e^{-1} + O(0), ax_e^k + bx_f^k+O(k+1)]&, ax_e^k +bx_f^k + O(k+1)] \\
&= [[c,a],a] x_e^{2k-1}+O(2k)
\in W
\end{align*}
But if $k>1$ then $2k-1 >k$. The lemma is proved.
\end{proof}

\begin{lemma}
\label{lemma:4}
For an arbitrary element $a \in S$ the subalgebra $W$ contains an element of the form $a x_f^{n-1} + O(n)$.
\end{lemma}

\begin{proof}
Let $a \in S$. Then $p_0(a) = p_1(a) =0$ and therefore $a x_e^{n-1} + p_2(a)
x^2 + p_3(a) x^3 \dots \in W$. Choose a minimal $k$ such that $p_k(a) \ne
0$. If $k \le n-1$, then $p_k(a)x^k + O(k+1) \in W$, which contradicts Lemma
\ref{lemma:3}.

If $k = n-1$, then $(a+ p_{n-1}(a))x_e^{n-1} + p_{n-1}(a) x_f^{n-1}+O(n) \in
W$. Again by Lemma \ref{lemma:2} $a+p_{n-1}(a) =0$, hence
$p_{n-1}(a)x_f^{n-1} + O(n) = -a x_f^{n-1}+O(n) \in W$. If $ k \ge n$, then
$a x_e^{n-1}+O(n) \in W$, which again contradicts Lemma \ref{lemma:3}.

Finally, if for any $k \ge 2$ we have $p_k(a) = 0$, then $a x_e^{n-1}\in W$,
which contradicts Lemma \ref{lemma:3}. The lemma is proved.
\end{proof}

\begin{lemma}
\label{lemma:5}
$S^\perp = \{a \in \g :\ \mathrm{\ there \ exists \ } b \in \g \mathrm{ \ such \ that \ } af + be + O(1) \in W\} = \{a \in \g :\ a + O(1) \in W \}$.
\end{lemma}

\begin{proof} Let $X_1 =\{a \in \g:\ \mathrm{\ there \ exists \ } b \in \g \mathrm{ \ such \ that \ } af + be + O(1) \in W\}, X_2 =\{a \in \g :\ a + O(1) \in W \}$. Obviously , $X_2 \subseteq X_1$. Since $S^\perp = \mathrm{im}\ p_0 + \mathrm{im} \ q_0$ it follows that $S^\perp \subseteq X_2$. Now, let $a \in X_1, af + be+ O(1) \in W$. Choose an arbitrary element $c \in S$. By Lemma \ref{lemma:4} there exists an element $c X_f^{n-1} + O(n)$ in $W$. Hence $(af + be + O(1) | cx_f^{n-1} + O(n)) = K(a,c) = 0, a \in S^\perp$. The lemma is proved. \end{proof}

\begin{lemma}
\label{lemma:6}
For an arbitrary element $a \in \g$ there exists an element of the form $ae + O(1)$ in $W$.
\end{lemma}

\begin{proof}
For an arbitrary elements $a \in \g$ there exist elements $h_i(a) \in \g, i
\ge 0$ such that $ae+ \sum_{i\ge0}h_i(a) x^i \in W$. Then $h_0(a)f +
(a+h_0(a))e + O(1) \in W$. By Lemma \ref{lemma:4} there exists an element of
the type $h_0(a) + O(1)$ in $W$. Now, $(ae+ \sum_{i \ge 0} h_i(a)x^i)-
(h_0(a) +O(1)) = ae + O(1) \in W$. The lemma is proved.
\end{proof}

\begin{lemma}
\label{lemma:7}
If $ax_e+bx_f + O(2) \in W, a,b \in \g$, then $a = 0$.
\end{lemma}

\begin{proof} Let $I = \{a \in \g:\ \mathrm{ \ there \ exists \ } b \in \g \mathrm{ \ such \ that } \ a x_e + b x_f + O(2) \in W\}, a \in I$. By Lemma \ref{lemma:6} for an arbitrary element $c \in \g$ there exists and element $ce+O(1) \in W$. Now $[ax_e + bx_f +O(2), ce  + O(1)] = [a,c]x_e + O(2) \in W$. Hence $I$ is an ideal in $\g$, hence $I = (0)$ or $I = \g$. Choose another pair of elements  $a' \in I, c' \in \g$. By the above there exists an element of the type $[a',c'] + O(2)$ in $W$. Hence \[ [ [a,c]x_e + O(2), [a',c']x_e + O(2)] = [[a,c],[a',c']]x_e^2 + O(3) \in W.\]

By Lemma \ref{lemma:2} $[[a,c],[a',c']] =0$ and thus $[[I,\g], [I,\g]] =(0)$
which implies that $I = (0)$. The lemma is proved.
\end{proof}

\begin{lemma}
\label{lemma:8}
For an arbitrary element $a \in \g$ there exists an element of the form $a x_f + O(2)$ in $W$.
\end{lemma}

\begin{proof} If $a \in \ker p_0$, then $p_1(a)x+ax_e^{n-1}+p_2(a)x^2+\dots = p_1(a)x+O(2) \in W$. Hence, by Lemma \ref{lemma:7} $p_1(a) =0$. We proved that $\ker p_0 \subseteq \ker p_1$, which implies that $\mathrm{im}\ p_0 = (\ker p_0)^{\perp} \supseteq (\ker p_1)^\perp = \mathrm{im} \ q_0$.

Hence for an arbitrary element $ a \in \g$ there exists $b \in \g$ such that
$p_0(b) = q_0(a)$. We have also \[q_0(a) + (q_1(a) +a)x_e + q_1(a)x_f + O(2)
\in W,\quad p_0(b) + p_1(b)x +O(2) \in W.\]

Subtracting these two inclusions we get \[ (q_1(a) + a-p_1(b)) x_e + (q_1(a)- p_1(b)) x_f + O(2) \in W.\]

By Lemma \ref{lemma:7} $q_1(a) - p_1(b) = -a$, hence $-a x_f +O(2) \in W$, which proves the lemma.
\end{proof}

\begin{proposition}
\label{prop:1}
If $\g \otimes A(n,\alpha)$ has a Lagrangian subalgebra $W$ such that
$\g \otimes A(n,\alpha) = \g \otimes F[[x]]+W$ is a direct sum, then $n \le 2$.
\end{proposition}

\begin{proof}
Assume that $n>2$. Since $\g = [\g,\g]$, commuting elements of the type
$ax_f + O(2) \in W$ we conclude that for an arbitrary element $a \in \g$
there exists an element of the type $a x_f^{n-2}+O(n-1) \in W$. Let $b \in
\g$, $bx_f + O(2)\in W$. Then \[(ax_f^{n-2}+O(n-1) | b x_f + O(2) ) = K(a,b)
=0.\] Hence $K(\g,\g) = (0)$, a contradiction. The proposition is proved.
\end{proof}

Now our aim is to prove the following.

\begin{proposition}
\label{prop:2}
If $A(2, \alpha)$ is a trace extension of $F[[x]]$ and $\g \otimes A(2,\alpha)$ has a Lagrangian subalgebra $W$, $\g \otimes A(2,\alpha) = \g \otimes F[[x]] +W$, then $\alpha_0 = 0$.
\end{proposition}

Recall that $A(2,\alpha) = F((x_e)) + Ff +F x_f$. Hence for an arbitrary element $a \in \g$ there exist elements $p_i(a), q_i(a) \in \g$, $i \ge 0$ such that
\begin{gather*}
ax_e  + p_0(a) + p_1(a)x + p_2(a)x^2 + \dots \in W \\
ae + q_0(a) + q_1(a)x + \dots \in W.
\end{gather*}
Denote $S = \ker p_0 \cap \ker q_0$.  Just as in the proof of Proposition
\ref{prop:1}
$$
\left(a x_e + \sum_{i \ge 0} p_i(a) x^i |\ b x_e + \sum_{i \ge 0} p_i(b)x^i\right)=0
$$
implies  $K(p_0(a),b) + K(a,p_0(b))=0$ and
$$
\left(ax_e+\sum_{i \ge 0} p_i(a) ^i |\ be + \sum_{i \ge 0} q_i(b) x^i\right) =0
$$
implies that
\begin{align*}
K(a,b) + K(a,q_0(b)) + \alpha_0K (p_0(a),b) + K (p_1(a), b)
&= K(a, (1+q_0)(b))\\
&\quad + K((\alpha_0 p_0+p_1)(a),b)\\
&=0.
\end{align*}
Hence,
$$
p_0^T = -p_0,\quad (\alpha_0p_0 + p_1)^T = -(1+q_0).
$$
This implies $S^\perp = \mathrm{im} \ p_0 + \mathrm{im} \ (1+q_0)$.

\begin{lemma}
\label{lemma:9}
$S = \{a \in \g:\ a x_e + O(2) \in W\}$.
\end{lemma}

\begin{proof}
The inclusion of $S$ in the right hand side is obvious. Now let $z = a x_e + O(2) \in W$. We have
$$
\g \otimes A(2,\alpha)
= \sum_{k=1}^\infty \g x_e^{-k} + \g e + \g x_e + \g\otimes F[[x]].
$$
Comparing degrees, we see that
$z \in \g e + \g x_e + \g \otimes F[[x]]$.
Since $z \in W$ it follows that
\[
z = (a_0 e + q_0(a_0) + q_1(a_0)x+\dots) + (a_1 x_e + p_0(a_1)+ p_1(a_1)x+\dots),
\]
where $a_0, a_1 \in \g$. Since $z \in x_e + O(2)$ it follows that $a_0e +
q_0(a_0) + p_0(a_1) =0$, $q_1(a_0)x+a_1x_e+p_1(a_1)x \in \g x_e$, which
implies $a_0 =0$, $p_0(a_1) =0$, $p_1(a_1) =0$. Therefore $a \in S$. The
lemma is proved.
\end{proof}

\begin{lemma}
\label{lemma:10}
$S^\perp = \{ a \in \g: \mathrm{ \ there \ exists \ } b \in \g \mathrm{, \ such \ that \ } ae + bf + O(1) \in W\}$.
\end{lemma}

\begin{proof} It is clear that $S^\perp = \mathrm{im} \ p_0 + \mathrm{im} (1 + q_0)$ lies in the right hand side.
Now suppose that $a \in \g, b \in \g$, and $ae + bf + O(1) \in W$. Let $c \in S$. Then by Lemma \ref{lemma:9} we have $cx_e + O(2) \in W$. Hence,
$$
(ae + bf + O(1) |\ c x_e + O(2)) = K(a,c) = 0,
$$
which implies the other inclusion. The lemma is proved.
\end{proof}

\begin{lemma}
\label{lemma:11}
$[S,S]=(0)$.
\end{lemma}

\begin{proof}
Let $a,b \in S$. Then there exist elements $a x_e + O(2)$,
$bx_e + O(2)$ lying in $W$. Hence,
\[ax_e + O(2),
bx_e +O(2)] = [a,b]x^2 + O(3) \in W \cap \g\otimes F[[x]] = (0), \] which
implies $[a,b] =0$. The lemma is proved.
\end{proof}

\begin{lemma}
\label{lemma:12}
$[\g,S] \subset \{ a \in \g:\  ae + O(1)\in W\} \subset S^\perp.$
\end{lemma}

\begin{proof} Choose an element $a \in S$, $a x_e + O(2) \in W$. For an arbitrary element $b \in \g$ there exists an element of the type $b x_e^{-1} +O(0)$ in $W$. Hence,
$$
[a x_e + O(2), bx_e^{-1}+O(0)] = [a,b]e + O(1) \in W.
$$
By Lemma \ref{lemma:10} we see that $[a,b] \in S^\perp$. The lemma is
proved. \end{proof}

\begin{lemma}
\label{lemma:13}
$[S^\perp, S] \subseteq S$.
\end{lemma}

\begin{proof}
Let $a \in S$, $a x_e + O(2) \in W$ and $b \in S^\perp$, $be+ cf+ O(1) \in W$ for some element $c \in \g$. Then
$$
[ax_e + O(2), be + cf + O(1)] = [a,b]x_e +O(2) \in W
$$
which implies the claim. The lemma is proved.
\end{proof}

Recall that an abelian subalgebra $I$ of a Lie algebra $L$ is called an
inner ideal if $[[L,I],I] \subseteq I$. Lemmas \ref{lemma:10},
\ref{lemma:11}, and \ref{lemma:12} imply that $S$ is an inner ideal of the
Lie algebra $\g$.

\begin{lemma}
\label{lemma:14}
If $S \ne (0)$ then $\alpha_0 =0$.
\end{lemma}

\begin{proof} If $S \ne (0)$ then $S$ contains a minimal nonzero inner ideal of $\g$. In \cite{B} it is proved that minimal inner ideals in a semisimple finite dimensional Lie algebra over an algebraically closed field of zero characteristic are one dimensional. Hence $S \ni a \ne 0, [a,[a,\g]] = Fa$. By the Morozov Lemma (see \cite{J}) there exists an element $b \in \g$ such that $a, h =[a,b], b$ is an $\mathfrak{sl}_2$-triple, that is $[h,a]=2a, [h,b] = -2b$. Since $h \in [a,\g] \in [S,\g]$ by Lemma \ref{lemma:12} it follows that there exists an element $he+ O(1)$ in $W$. Choose elements $c,d \in \g$ such that $he+ cx_e + dx_f +O(2) \in W$. Since $a \in S$ there exists an element $a x_e + O(2)$ in $W$. Choose an element $v \in \g$ such that $ax_e + vx^2 + O(3) \in W$. We have \[ [he + cx_e + dx_f + O(2), ax_e + vx^2 + O(3)] = 2ax_e + ([h,y]+[c,a])x^2 + O(3) \in W.\] Now \[ (2ax_e + ([h,y]+[c,a])x^2 + O(3)) -2(ax_e +vx^2+O(3)) \in W \cap \g \otimes F[[x]] = (0).\]

Hence $[h,y]+[c,a] = 2v$. This implies $K(b,[h,y]) + K(b,[c,a]) = 2K(b,y)$. But $K(b,[h,y])=K([b,h],y) = 2K(b,y)$. Hence $K(b,[c,a]) = K(h,c)=0$.

On the other hand since the element $he + cx_e + dx_f + O(2)$ lies in $W$,
we have $(he+ cx_e + d x_f + O(2) | he + c x_e + d x_f + O(2)) = \alpha_0
K(h,h) + 2 K(h,c) = 0$. Since $K(h,h) \ne 0$ it follows that $\alpha_0 =0$.
The lemma is proved.
\end{proof}

\begin{lemma}
\label{lemma:15}
$\ker p_0 \ne (0)$.
\end{lemma}

\begin{proof} If $\ker p_0 = (0)$, then $p_0$ is invertible. Denote $f = p_0^{-1}$. Since $W$ is a subalgebra we have
\begin{align*}
\left[a x_e + \sum_i p_i(a) x^i, bx_e + \sum_i p_i(b) x^i \right]
&= ([a,p_0(b)]+[p_0(a),b])x_e \\
&\quad + [p_0(a),p_0(b)]+ \sum_{i \ge 1}(\dots) x^i
\end{align*}
Therefore $p_0([p_0(a),b]+[a,p_0(b)]) = [p_0(a),p_0(b)]$ or, equivalently,
$f([a,b]) = [f(a),b] + [a,f(b)]$. Hence, $f$ is a derivation, which
contradicts it being invertible. The lemma is proved. \end{proof}

Now our aim is to prove that $\alpha_0 =0$. We will assume therefore the contrary. In particular $S=(0)$ by Lemma \ref{lemma:14}.

\begin{lemma}
\label{lemma:16}
For an arbitrary element $a \in \ker p_0$ there exists an element of the form $ax_f + O(2)$ in $W$.
\end{lemma}

\begin{proof} Let $a \in \ker p_0$. Then
$$
ax_e + p_1(a) + O(2) = (a + p_1(a)) x_e + p_1(a) x_f + O(2) \in W.
$$
For an arbitrary element $b \in \g$ choose $b x_e^{-1} + O (0) \in W$. Then \begin{align*} [ [ bx_e^{-1} + O(0), (a+p_1(a))x_e  +  & p_1  (a)x_f +   O(2)], (a +p_1(a))x_e + p_1(a) x_f+O(2)]  \\ & = [[b, a+p_1(a)], a+p_1(a)]x_e +O(2) \in W. \end{align*}
Hence, $[[\g, a+p_1(a)], a+p_1(a)] \subseteq S = (0)$, which implies that $a+p_1(a) =0$. Now,
$$
(a+p_1(a))x_e + p_1(a)x_f + O(2) = -ax_f +O(2) \in W.
$$
The lemma is proved.
\end{proof}

\begin{lemma}
\label{lemma:17} For an arbitrary element $a \in \g$ there exists a unique
element $ b \in \g$ such that $ae + bx_e + O(2) \in W$.
\end{lemma}

\begin{proof}
Let $T= \{a \in \g: ae + O(1) \in W\}$. Since $S = (0)$, $S^\perp = \g$. By Lemma \ref{lemma:10} it follows that for an arbitrary element $c \in \g$ there exists $d \in \g$ such that $ce+ df + O(1) \in W$. Let $a \in T$. Then
$$
[ae+O(1), ce+df+O(1)] = [a,c]e+O(1)\in W.
$$
Hence $[T,\g] \subset T$. Hence $T =\g$ or $T = (0)$. Suppose that $T= (0)$. For an arbitrary element $a \in \g$ we have $ae + q_0(a) + O(1) \in W$. Hence $ \ker q_0 \subset T = (0)$, so $q_0$ is invertible. This implies that for an arbitrary element $b \in \g = \mathrm{im}\ q_0$ there exists an element $c \in \g$ such that $bf + ce+ O(1) \in W.$

If $a \in \ker p_0$ then by Lemma \ref{lemma:16} $a x_f + O(2) \in W$. Now  $(a x_f + O(2) | bf + ce + O(1)) = K(a,b) = 0$, which implies $ K(\ker p_0, \g)=(0)$. This contradicts Lemma \ref{lemma:15}. We have shown that $T = \g$. Now it is easy to see that for an arbitrary element $a \in [T,T] =\g$ there exists an element $b \in \g$, such that $ae+bx_e+O(2) \in W$.

Now, if $ae + bx_e + O(2)$ and $ae+b'x_e + O(2)$ both belong to $W$, then
$(b-b')x_e + O(2)\in W$, which implies $b-b' \in S = (0)$. The lemma is
proved.
\end{proof}

Now we are ready to finish the proof of Proposition \ref{prop:2}.
\begin{proof}
We will show that the assumption $\alpha_0 \ne 0$ leads to a contradiction. By Lemma \ref{lemma:17} there exists a map $h: \g \rightarrow \g$ such that for an arbitrary $a \in \g$ the element $ae + bx_e + O(2)$ lies in $W$ if and only if $b = h(a)$. Commuting two such elements, we get
\[
[ae + h(a)x_e+O(2), be + h(b) x_e + O(2)]
= [a,b]e + ([h(a),b] + [a,h(b)])x_e + O(2) \in W.
\]
Hence, $h([a,b]) = [h(a),b] + [a,h(b)]$,  i.e. $h$ is a derivation of $\g$. On the other hand, because of the isotropy of $W$ we have
$$
(ae +h(a) x_e + O(2) | be + h(b)x_e + O(2))
= \alpha_0 K(a,b) + K (h(a),b) + K(a,h(b)) =0.
$$
Since $h$ is a derivation it follows that $K(h(a),b)  + K(a,h(b)) =0$. Thus if $\alpha_0 \ne 0$ then $K(a,b) =0$ for arbitrary $a,b \in \g$, which is a contradiction. Proposition \ref{prop:2} is proved.
\end{proof}

\subsection*{End of the proof of Theorem \ref{thm:1}}

\begin{proof}

If $A = A(2,\alpha)$, $\alpha_0 =0$ then it is not difficult to find coefficients $\xi_i \in F, i \ge 1$ such that for $y_e  = x_e + \sum_{i \ge 1}\xi_i x_e^{i+1}$ we have $t(y_e^{-k}) =0$ for $k=1,2, \dots$
Let
$$
x_e = y_e + \sum_{i \ge 1} \eta_i y_e^{i+1}
\quad{\rm and} \quad
\phi(x) = x + \sum_{i\ge1} \eta_i x^{i+1}.
$$
Then $A(2,\alpha) \cong A(2,0)^{\phi(0)}$.
If $n=1$ then
$$
A(1,\alpha) = F((x_e)) \oplus Ff, \quad t(e)=1,\quad t(f)=-1.
$$
Again there exists a series
\[y_e = x_e + \sum_{i \ge 1} \xi_i x_e^{i+1},\quad
x_e = y+e + \sum_{i \ge 1} \eta_i y_e^{i+1}\]
such that $t (y_e^{-k}) =0$ for $k=1,2, \dots$.
Then $A(1,\alpha) \cong A(1,0)^{(\phi)}$ for
$\phi(x) = x + \sum_{i \ge 1} \eta_i x^{i+1}$.

In the case of $n=0$ we have $A(0,\alpha) = F((x))$, $t(x^{-1}]=1$ and there exist series $y = x + \sum_{i \ge 1} \xi_i x^{i+1}$, $x = y+ \sum_{i\ge1}\eta_i y^{i+1}$ such that $t(y^{-k})=0$ for $k=2,3,\dots$. Then
$$
A(0,\alpha) \cong A(0,0)^{(\phi)},\quad \phi(x) = x + \sum_{i \ge 1}\eta_i x^{i+1}.
$$
This finishes the proof of Theorem \ref{thm:1}.
\end{proof}

\begin{remark}
Let us consider the group of $F[[x]]$-linear automorphisms of $\g[[x]]$, which we will denote by $Aut_{F[[x]]}(\g[[x]])$.
Each automorphism $U\in Aut_{F[[x]]}(\g[[x]])$ can be extended to an automorphism $\tilde{U}$ of $A(n,\alpha)$.
The automorphism $\tilde{U}$ preserves the sequence $\{\alpha_i\}$.
\end{remark}

\section{Multivariable case}

Let us now consider the algebra of the series $F[[X]] = F[[x_1, \dots, x_n]]$ of $n \ge 2$ variables. Let $F[[X]] \subseteq A$, $t:A\rightarrow F$, be a trace extension, that is, the bilinear form $(a|b)=t(ab)$ on $A$ is nondegenerate,
$$
F[[X]]^\perp = \{a \in A:\ (F[[X]] |a) = (0) \} = F[[X]]
$$
and for an arbitrary distribution $f:F[[X]] \rightarrow F$ there exists an element $a \in A$ such that $f(b) = (b|a)$ for an arbitrary $b \in F[[X]]$.

\begin{theorem}
\label{thm:4}
The trace extension $(A,t)$ is isomorphic to a trivial extension.
\end{theorem}

\begin{proof} For a multi-index $\alpha = (\alpha_1, \dots, \alpha_n) \in \mathbb{Z}_{\ge0}^n$ denote $x^\alpha = x_1^{\alpha_1}\dots x_n^{\alpha_n}$. Choose elements $b_\alpha \in , \alpha \in \mathbb{Z}_{\ge 0}^n$ such that $t(b_\alpha x_\beta) = \delta_{\alpha,\beta}$, the Kronecker symbol. Our aim is to construct elements $a_\alpha \in b_\alpha +A, \alpha \in \mathbb{Z}_{\ge 0}^n$ such that \[a_\alpha x^\beta = \begin{cases}a_{\alpha-\beta} & \mathrm{if} \ \alpha \ge \beta \\ 0 & \mathrm{otherwise.} \end{cases}\]
We will proceed by induction on $|\alpha| = \alpha_1 + \dots + \alpha_n$. To start we assume $a_\alpha=0$ for all $\alpha \in \mathbb{Z}^n \setminus \mathbb{Z}_{\ge 0}^n$. Let $\epsilon(i) = (0, \dots, 1, \dots, 0)$ (with 1 on the $i$th place).

Choose $\alpha \in \mathbb{Z}_{\ge 0}^n$ and suppose that the elements $a_{\alpha - \epsilon(i)}$ have been chosen. For an arbitrary $i, 1 \le i \le n$ we have $(b_\alpha x_i -a_{\alpha-\epsilon(i)}|F[[X]]) =(0)$. This implies that $b_\alpha x_i - a_{\alpha - \epsilon(i)} = \g_{\alpha, i} \in F[[X]]$. Furthermore, $b_\alpha x_i x_j - a_{\alpha - \epsilon(i) - \epsilon(j)} = \g_{\alpha,i}x_j = \g_{\alpha,j}x_i$ for $i \ne j$. Hence there exists an element $p_\alpha \in F[[X]]$ such that $\g_{\alpha,i} = p_\alpha x_i$. We then set $a_\alpha = b_\alpha - p_\alpha$. It is easy to see that the elements $a_\alpha, \alpha \in \mathbb{Z}_{\ge 0}^n$ satisfy the above assumptions. Choose $\alpha, \beta \in \mathbb{Z}_{\ge 0}^n$ and choose $\gamma \in \mathbb{Z}_{\ge 0}^n$ such that $\gamma > \alpha$. Then $a_\beta = a_{\beta+\gamma} x^\gamma$ and $a_\alpha a_\beta = a_\alpha x^\gamma a_{\beta+\gamma} = 0$.
It is easy to see that $A = F[[X]] +  \sum F a_\alpha$ is the trivial trace extension of $F[[X]]$. Theorem \ref{thm:4} is proved.
\end{proof}

\section{Classification over polynomials}

The aim of this part is to use Theorem \ref{thm:1} to classify classical
doubles over polynomials. In this section we assume that $F=\mathbb{C}$.



\begin{lemma}
\label{lemma:1a}
Let $\delta:\g[x] \rightarrow \g[x] \otimes \g[y] = (\g \otimes \g) [x,y]$ be a Lie bialgebra structure on $\g[x]$. Then it can be extended to $\overline{\delta}:\g[[x]] \rightarrow (\g \otimes \g)[[x,y]]$.
\end{lemma}

\begin{proof} Since $\g = [\g,\g]$, we can write $ax^2 = [bx,cx]$ where $a,b,z \in \g$.  Hence
\begin{align*}
\delta(ax^2)
&= \delta([bx,cx])\\
&=[\delta(bx),cx\otimes 1 +1 \otimes cy]
+ [bx \otimes1 +1 \otimes by, d(cx)] \in (x,y) \cdot (\g \otimes \g)[x,y].
\end{align*}
A simple induction shows that
$$
\delta(ax^n) \in (x,y)^{n-1}(\g \otimes \g)[x,y].
$$
Therefore $\overline{\delta}(\sum_{n=0}^\infty a_n x^n)$ can be defined as $\sum_{n=0}^\infty \delta(a_n x^n)$.
\end{proof}

It follows from Theorem \ref{thm:1} that the classical double $D(\g[[x]],
\overline{\delta}) \cong \g \otimes A(n,\alpha)$ where $n =0,1,2$. Our plan
is to use the condition $\overline{\delta}(\g[x]) \subseteq (\g \otimes
\g)[x,y]$ to get a description of all possible sequences $\alpha$. Since the
proofs in all three cases are similar, we will provide detailed proofs in
case $n=0$.

In this case $\g \otimes A(0,\alpha) \cong \g((x))$ and the trace $t:
\mathbb{C}((x)) \rightarrow \mathbb{C}$ is given by the formula $t(x^n) =0$,
$n \ge 0$ and $t(x^{-1}) =1$, $t(x^{-k}) = a_{k-1}$ for $k \ge 2$. Denote $1
+ \sum_{i \ge 1}a_i x^i$ by $a(x)$. Then, evidently the canonical form in
the double $A(0,\alpha)$ is given by the formula $(f_1(x)|f_2(x))=
\mathrm{Res}_{x=0}(K'(f_1,f_2))\cdot a(x))$, where $K'$ is the Killing form
of the Lie algebra $\g((x))$ over $\mathbb{C}((x))$.

\begin{lemma}
\label{lemma:2a} Let $W \subset \g \otimes A(0,\alpha)$ be a subspace. Let
$W^{\perp_0}$ be the orthogonal complement of $W$ in $\g \otimes A(0,\alpha)
= \g((x))$ with respect to $(f_1|f_2)_0 = \mathrm{Res}_{x=0} K'(f_1,f_2)$.
Then the orthogonal complement $W^\perp$ of $W$ with respect to the form
$(f_1|f_2) = \mathrm{Res}_{x=0}(K'(f_1,f_2)\cdot a(x))$ is
$\frac{1}{a(x)}W^{\perp_0}$.
\end{lemma}

\begin{theorem}
\label{thm:1a} Let $\g \otimes A(0,\alpha) =g [[x]] \oplus W$, the $W$ is a
Lagrangian subalgebra of $A(0,\alpha)$ corresponding to $\overline{\delta}$.
Then $W$ is bounded, i.e., there exists $N$ such that $W \subseteq x^N
\g[x^{-1}]$.
\end{theorem}

\begin{proof} Let $\{E_k\}$ be an orthonormal basis of $\g$ with respect to the Killing form.
Since $\g \otimes A(0,\alpha) = \g[[x]] \oplus x^{-1}\g[x^{-1}]$,
there exists a basis of $W$, which consists of the following elements:
\[
E_{k,n}=E_k x^{-n}+p_{k,n}(x), n =-1,-2, \dots,
\]
where $p_{k,n} \in \g[[x]]$. By Lemma \ref{lemma:2a} the dual basis in $\g[[x]]$ with respect to the form $(f_1,f_2) =\mathrm{Res}_{x=0}(K'(f_1,f_2) \cdot a(x))$ consists of the elements $\frac{1}{\alpha(x)}E_k x^{n-1}\in \g[[x]]$. Then
\[
\overline{\delta}(f(x)) = [f(x) \otimes 1 + 1 \otimes f(y), r(x,y)] = \frac{1}{a(x)}
 \left ( \frac{\Omega}{y-x} + R(x,y) \right ),
 \]
where $\Omega = \sum E_k \otimes E_k\in \g \otimes g$ is invariant and $R(x,y) \in (\g \otimes \g)[[x,y]]$. Let
$$
\frac{1}{a(x)} =1 + b_1 x + \dots + b_n x^n \cdots
$$
Then
\[
\frac{1}{a(x)} R(x,y)
= \sum_{k,n} E_k x^{n-1} (1 + b_1x + \dots + b_k x^k \dots)
\otimes p_{k,n}(y)
\]
is well defined because  \[ \frac{1}{a(x)} R(x,y)  = \sum_{n,k} E_k x^n
\otimes \bigg \{ \mathrm{ \ finite \ sum \ of \ } \ p_{\alpha, k}(y) \bigg \} \]
Let us rewrite $\frac{1}{a(x)} R(x,y)$ as $\sum c_{ij} x^i y^j \in (\g \otimes \g) [[x,y]]$ and let us compute $\overline{\delta}(E)$ for $E \in \g$. We get
$\overline{\delta}(E) = \sum[c_{ij},E\otimes 1 + 1 \otimes E] x^i y^j$.
Since $\overline{\delta}(E) \in (\g\otimes \g)[x,y]$ we deduce that there exists a number $N$ such that $[c_{ij}, E\otimes 1 + 1 \otimes E] =0$ for all $i,j \ge N$.
It follows that $\frac{1}{a(x)} R(x,y) = P(x,y) + p(x,y)\Omega$, where $P \in (\g \otimes \g)[x,y]$ and $p(x,y) \in \mathbb{C}[[x,y]]$.
Now, we can rewrite
\[
r(x,y) = \frac{\Omega}{a(x)(y-x)} + \frac{(y-x)p(x,y)\Omega}{y-x} + P(x,y) = \frac{A(x,y)\Omega}{y-x} + P(x,y).
\]
Here $A(x,y) \in \mathbb{C}[[x,y]]$.

Now, let us compute for $E \in \g$,
$$
\overline{\delta}(Ex)
= A(x,y)[\Omega, E\otimes1]+ [P(x,y),Ex\otimes 1 + 1 \otimes Ey].
$$
Since the second summand is polynomial and $\overline{\delta}(Ex)$ is polynomial, we have proved that $A(x,y)$ is polynomial.

Let us notice that
$$
\frac{x^iy^j}{y-x} = \sum_{k \ge 0} x^{k+i}y^{-k-1+j}
$$
and therefore we can rewrite $r(x,y) = \sum_{k,n}E_k x^n \otimes
s_{k,n}(y)$, where $s_{k,n}(y) \in y^N \g[y^{-1}]$ with $N$ being maximal of
the degrees of $A(x,y)$ and $P(x,y)$ in $y$. Let $\{E_k x^n\}$ and
$E_{k,n}(y)$ be the dual bases of $\g[[x]]$ and $W$ with respect to
$(f_1,f_2)=\mathrm{Res}_{x=0}(K'(f_1,f_2)\cdot a(x))$. Then $r(x,y) = \sum
E_k x^n \otimes E_{k,n}(y)$ and we see that $s_{k,n}(y) = E_{k,n}(y)$ and $W
\subset y^N \g[y^{-1}]$. The theorem is proved.
\end{proof}

\subsection*{Maximal orders in loop algebras}

In what follows we need the so-called orders in $\g((x^{-1}))$. Proofs of the results below can be found in \cite{S2,S3}.

\begin{definition}
Let $W \subseteq \g((x^{-1}))$ be a Lie subalgebra. We say that $W$ is an order if there exists $n,k$ such that $x^{-n} \g[[x^{-1}]] \subseteq W \subseteq x^k \g[[x^{-1}]]$.
\end{definition}

\begin{example} $W_0=\g[[x^{-1}]]$.
\end{example}

Let us consider the group $\mathrm{Aut}_{\mathbb{C}[x]}(\g[x])$. Clearly, there exists a natural embedding $i:\mathrm{Aut}_{\mathbb{C}[x]}(\g[x]) \hookrightarrow \mathrm{Aut}_{\mathbb{C}((x^{-1}))}(\g((x^{-1}))$. If $\sigma(x) \in \mathrm{Aut}_{\mathbb{C}[x]}(\g[x])$, then abusing notations we denote $i(\sigma(x))$ by $\sigma(x)$.

\begin{definition} We say that two orders $W_1$ and $W_2$ are gauge equivalent if there exists $\sigma(x) \in \mathrm{Aut}_{\mathbb{C}[x]}(\g[x])$ such that $\sigma(x)W_1 = W_2$.
\end{definition}

Now, let us introduce some orders in $\g((x^{-1}))$, which will be denoted by $\mathbb{O}_h$.

Fix a Cartan subalgebra $\h \subset \g$. Let $R$ be the corresponding set of
roots and $\Gamma$ the set of simple roots. Denote by $\g_\alpha$ the root
space corresponding to $\alpha \in R$. Consider the valuation on
$\mathbb{C}((x^{-1}))$ defined by $v(\sum_{k \ge n}a_k x^{-k}) =n$,
$a_n\neq0$. For every root $\alpha \in R$ and every $h \in \h(\mathbb{R})$,
set
$$
M_\alpha(h):=\{ f \in \mathbb{C}((x^{-1})): v(f) \ge \alpha(h)\}.
$$
Consider $\mathbb{O}_h:= \h[[x^{-1}]]\oplus \left ( \bigoplus_{\alpha \in R}
M_{\alpha}(h)\g_\alpha \right)$. It is not difficult to see that
$\mathbb{O}_h$ is a Lie subalgebra and, moreover, an order in
$\g((x^{-1}))$.

Let us consider the following standard simplex
$$
\{h \in \eta(\mathbb{R}): \alpha(h) \ge 0
\quad \forall
\alpha \in \Gamma, \quad \alpha_\mathrm{max}(h) \leq 1 \}.
$$
The vertices of this simplex are 0 and $h_i$, where $h_i$ are uniquely defined by the condition $\alpha_j (h_i) = \delta_{ij}/k_i$. Here $\alpha_{\mathrm{max}} = \sum k_j \alpha_j$, $\alpha_j \in \Gamma$. Clearly, there exists a one-to-one correspondence between the vertices of the standard simplex and the vertices of the extended Dynkin diagram of $\g:\alpha_0:= \alpha_{\mathrm{max}} \leftrightarrow 0$, $ \alpha_i \leftrightarrow h_i$. In what follows we will denote $\mathbb{O}_{h_i}$ by $\mathbb{O}_{\alpha_i}$.

\begin{theorem}
\label{thm:2a}
Let $W \subset \g((x^{-1}))$ be an order such that $W + \g[x] = \g((x^{-1}))$. Then there exists $\sigma(x) \in \mathrm{Aut}_{\mathbb{C}[x]}(\g[x])$ such that $\sigma(x)W \subseteq \mathbb{O}_{\alpha_i}$ for some $\alpha_i \in \Gamma$. \end{theorem}

\begin{remark}
Generally speaking $\alpha_i$ is not unique, it may happen that there exists $\sigma_1(x) \in \mathrm{Aut}_{\mathbb{C}[x]}(\g[x])$ such that $\sigma_1(x)W \subseteq \mathbb{O}_{\alpha_j}$.
\end{remark}

\subsection*{Application of maximal orders to classification of Lie bialgebra structures on $\g[x]$}

We continue with the case $D(\g[[x]], \overline{\delta}) \cong \g\otimes
A(0,\alpha)$, where $\overline{\delta}$ is a natural extension of $\delta:
\g[x] \rightarrow (\g \otimes \g)[x,y]$. Let $\g\otimes A(0, \alpha) =
\g[[x]] \oplus W$ be the Manin triple corresponding to $\overline{\delta}$.
Then according to Theorem \ref{thm:1a}, $W \subseteq x^N \g[x^{-1}]$ for
some $N$. Since $W = W^\perp$, we see that $(x^N \g[x^{-1}])^\perp \subset W
\subset x^N\g[x^{-1}]$.

\begin{lemma}
\label{lemma:3a}
\[ (x^N \g[x^{-1}])^\perp = \frac{1}{a(x)} x^{-N-2}\g[x^{-1}] \]
\end{lemma}
\begin{proof}
Let us recall that for any $V \subset \g\otimes A(0,\alpha)$, $V^\perp =
\frac{1}{a(x)}V^{\perp_0}$. It remains to notice that $(x^N
\g[x^{-1}])^{\perp_0} = x^{-N-2}\g[x^{-1}]$.
\end{proof}

\begin{corollary}
\label{cor:1a}
$\frac{1}{a(x)}$ is a polynomial of the form $1+ b_1 x + \dots b_k x^k$, $k \le 2N +2$.
\end{corollary}

\begin{proof}
We have
$$
\frac{x^{-N-2}}{a(x)} \g [x^{-1}] \subset W^\perp = W \subset x^N \g[x^{-1}],
$$
which proves the corollary.
\end{proof}

In order to formulate the main result of this section let us make some remarks.

\begin{enumerate}[1.]
\itemsep-3pt
\item
$W \subset \g[x,x^{-1}] \subset \g ((x^{-1}))$.
\item
$W \cdot \mathbb{C}[[x^{-1}]]$ is an order in $\g((x^{-1}))$
because $W \cdot \mathbb{C}[[x^{-1}]] \subset x^N \g [[x^{-1}]]$
and
\begin{align*}
W \cdot \mathbb{C}[[x^{-1}]]
&\supset x^{-N-2}(1+b_1 x+ \dots+b_k x^k) \cdot \g[[x^{-1}]] \\
&\supset x^{-N-k-2}(b_k+b_{k-1}x^{k-1} + \dots + x^{-k}) \cdot \g[[x^{-1}]]\\
&= x^{-N-k-2} \g[[x^{-1}]]
\end{align*}
since $b_k +b_{k-1} x^{-1} + \dots + x^{-k}$
is a unit in $\mathbb{C}[[x^{-1}]]$.
\item
Therefore, $W \cdot \mathbb{C}[[x^{-1}]]$ is an order in $\g((x^{-1}))$.
\item
There exists $\sigma(x) \in \mathrm{Aut}_{\mathbb{C}[x]}(\g[x])$ such that $\sigma(x)(W\cdot \mathbb{C}[[x^{-1}]])\subset \mathbb{O}_{\alpha_i}$.
\end{enumerate}

\begin{theorem}
\label{theorem-ma}
Assume that $\sigma(x)(W \cdot \mathbb{C}[[x^{-1}]]) \subset \mathbb{O}_{\alpha_i}$.
Then
\begin{itemize}
\item
$\frac{1}{a(x)}$ is a polynomial of degree at most 2 if
$k_i =1$,
\item
$\frac{1}{a(x)}$ is a polynomial of degree at most 1 if $k_i >1$.
\end{itemize}
\end{theorem}

\begin{proof} Clearly $\sigma(x)W = W_1$ also defines a
Manin triple $\g\otimes A(0,\alpha) = \g[[x]] \oplus W_1$ and $W_1 \subset
\g[x,x^{-1}]$. It follows that $W_1 \subset \g[x,x^{-1}] \cap
\mathbb{O}_{\alpha_i}$. Let us describe $\g[x,x^{-1}] \cap
\mathbb{O}_{\alpha_i}$. For each $r, -k_i \le r \le k_i$, let $R_r$ denote
the set of all roots, which contain $\alpha_i$ with coefficient $r$. Let
$\g_0 = \h \oplus \sum_{\beta \in R_0} \g_\beta$ and $\g_r = \sum_{\beta \in
R_r} \g_\beta$. Then
\[\mathbb{O}_{\alpha_i} =  \sum_{r=1}^{k_i} x^{-1} \mathbb{C}[[x^{-1}]]\g_r + \sum_{r=1-k_i}^0 \mathbb{C}[[x^{-1}]]\g_r + x\mathbb{C}[[x^{-1}]]\g_{-k_i}\]
Hence
$$
\mathbb{O}_{\alpha_i} \cap \g[x,x^{-1}] = \sum_{r=1-k_i}^0 \mathbb{C}[[x^{-1}]]\g_r + \sum_{r=1-k_i}^0\mathbb{C}[x^{-1}]\g_r +x \mathbb{C}[x^{-1}]\g_{-k_i}.
$$
It is not hard to compute $(\mathbb{O}_{\alpha_i} \cap \g[x,x^{-1}] )^{\perp_0}$:
we will get \[x^{-3}\mathbb{C}[x^{-1}]\g_{k_i} + \sum_{r=0}^{k_i -1} x^{-2} + \sum_{r=-1}^{-k_i}\mathbb{C}[x^{-1}]\g_r.\]
By Lemma \ref{lemma:2a} we have
\begin{align*} \frac{1}{a(x)} & \left ( x^{-3} \mathbb{C}[x^{-1}] \g_{k_i} + \sum_{r=0}^{k_i -1} x^{-2} \mathbb{C}[x^{-1}] \g_r + \sum_{r=-1}^{-k_i} x^{-1} \mathbb{C}[x^{-1}] \g_r \right ) \subset \\ & W_1 \subset \left ( \sum_{r=1}^{k_i} x^{-1}\mathbb{C}[x^{-1}] \g_r + \sum_{r=1-k_i}^0 \mathbb{C}[x^{-1}]\g_r + x \mathbb{C}[x^{-1}] \g_{-k_i} \right ).
\end{align*}
We see that if $k_i = 1$, then
$$
(\mathbb{O}_{\alpha_i} \cap \g[x,x^{-1}])^{\perp_0}
= x^{-2} (\mathbb{O}_{\alpha_i} \cap \g[x,x^{-1}])
$$
and
$$
x^{-2}\frac{1}{a(x)} \left ( \mathbb{O}_{\alpha_i} \cap \g[x,x^{-1}]  \right ) \subseteq \mathbb{O}_{\alpha_i} \cap \g[x,x^{-1}]
$$
Then $\frac{1}{a(x)}$ is a polynomial of degree at most 2.

If $k_i >1$, we can consider the $\g_{k_i-1}$-component.
It is easy to see that $\frac{1}{a(x)}x^{-2}\mathbb{C}[x^{-1}]\g_{k_i -1} \subseteq x^{-1}\mathbb{C}[x^{-1}]\g_{k_i-1}$.
Therefore, $\deg \left ( \frac{1}{a(x)} \right) \le 1$.
The Theorem is proved.

\end{proof}

\begin{corollary}
\label{corollary-2} If
$D(\g[[x]],\overline{\delta})=\g\otimes{A}(0,\alpha)$, then
${D}(\g[x],\delta)\cong \g[x,x^{-1}]$ as Lie algebras. The canonical form on
$\g[x,x^{-1}]$ is given by the formula
$(f_{1}(x),f_{2}(x))=Res_{x=0}(K'(f_{1},f_{2})a(x))$. Up to automorphism
$\gamma:\mathbb{C}[x,x^{-1}]\longrightarrow \mathbb{C}[x,x^{-1}]$ given by
$\gamma(x)=cx$, $c\in\mathbb{C} $, we have the following possibilities for
$a(x)$:
$$a(x)=1\ (case\ A1)$$
$$a(x)=\frac{1}{1-x}\ (case\ A2)$$
$$a(x)=\frac{1}{(1-x)^{2}}\ (case\ A3)$$
$$a(x)=\frac{1}{(1-m_{1}x)(1-m_{2}x)}\ (case\ A4)$$
In the last case $m_{1}\neq m_{2}$, $m_{1}m_{2}\neq 0$, and the
corresponding canonical forms are parameterized by $\frac{m_{1}}{m_{2}}\in
(\mathbb{CP}^{1}\setminus \{0,1,\infty\})/\mathbb{Z}_2$.
\end{corollary}

The case $A1$ is well known and $W$ can be chosen as $W_{1}=x^{-1}\g[x^{-1}]$. In the case $A2$, one can easily check that $W$ can be chosen as
$$
W_{2}={\rm Span}\{(x^{-1}-\frac{1}{2})h_{i},(x^{-1}-1)e_{-\alpha},x^{-1}e_{\alpha},(1-x^{-1})x^{-1}\g[x^{-1}]\}.
$$
In the case $A3$ , one can check that the complementary $W$ can be chosen as
$$
W_{3}={\rm Span}\{(x^{-1}-1)e_{\pm \alpha},(x^{-1}-1)h_{i},(1-x)^{2}x^{-2}\g[x^{-1}]\}.
$$
Finally, for $A4$ one can check that
\begin{align*}
W_{4,m_{1},m_{2}}&=
{\rm Span}\{(x^{-1}-m_{1})e_{\alpha},(x^{-1}-m_{2})e_{-\alpha},
(x^{-1}-\frac{m_{1}+m_{2}}{2})h_{i},\\
&\qquad\qquad(1-m_{1}x)(1-m_{2}x)x^{-2}\g[x^{-1}]\}
\end{align*}
is a Lagrangian subalgebra complementary to $\g[x]$.

Let $\sigma$ be the Cartan involution of $\g[x,x^{-1}]$, $\sigma(e_{\alpha})=e_{-\alpha}$ for simple roots $\alpha_{i}$. Then $\sigma(\g[x])=\g[x]$ and
$\sigma(W_{4,m_{1},m_{2}})=W_{4,m_{2},m_{1}}$,
that implies that the forms corresponding to $m_{1}/m_{2}$ and
$m_{2}/m_{1}$ provide isomorphic Lie bialgebras on $\g[x]$.

Therefore we conclude that the variety of Lie bialgebras of type $A4$ is
isomorphic to $(\mathbb{CP}^{1}\setminus \{0,1,\infty\})/\mathbb{Z}_{2}$.

\begin{corollary}
\label{corollary-3}
The corresponding $r$-matrices are:
$$r_{A1}(x,y)=\frac{\Omega}{y-x}$$
$$r_{A2}(x,y)=\frac{1-x}{y-x}\Omega-r_{DJ}$$
$$r_{A3}(x,y)=\frac{(x-1)(y-1)}{y-x}\Omega$$
$$r_{A4,m_{1},m_{2}}(x,y)=\frac{1-(m_{1}+m_{2})u+m_{1}m_{2}xy}{y-x}\Omega-r_{m_{1},m_{2}}$$
\end{corollary}

Here $r_{m_{1},m_{2}}=\sum_{\alpha>0}m_{1}e_{-\alpha}\otimes e_{\alpha}+m_{2}e_{\alpha}\otimes e_{-\alpha}+\frac{m_{1}+m_{2}}{4}h_{i}\otimes h_{i}$
and $r_{DJ}=\frac{1}{2}(r_{-1,1}+\Omega)$,
which is the 
Drinfeld-Jimbo $r$-matrix.

The dual bases are given below:

\noindent Case A1:
$$
\left\{e_{-\alpha}y^{-k-1},\ e_{\alpha}y^{-k-1},\frac{1}{2}h_{i}y^{-k-1}\right\},\quad
k=0,1,\ldots
$$
Case A2:
\begin{gather*}
\{(y^{-1}-1)e_{-\alpha},\ y^{-1}e_{\alpha},  \  \frac{1}{2}(y^{-1}-\frac{1}{2})h_{i},\
(y^{-1}-1)y^{-k}e_{-\alpha},\\
(y^{-1}-1)y^{-k}e_{\alpha},\ \frac{1}{2}(y^{-1}-1)y^{-k}h_{i}\},\quad
k=1,2,\ldots
\end{gather*}
Case A3:
$$
\{ (y^{-1}-1)e_{-\alpha}, \ (y^{-1}-1)e_{\alpha},\ \frac{1}{2}(y^{-1}-1)h_{i},\ (1-y)^{2}y^{-k-1}e_{-\alpha},$$
$$
(1-y)^{2}y^{-k-1}e_{\alpha}, \ \frac{1}{2}(1-y)^{2}y^{-k-1}h_{i}\},\quad
k=1,2,\ldots
$$
Case A4:
$$
\{(y^{-1}-m_{2})e_{-\alpha},\ (y^{-1}-m_{1})e_{\alpha},\ \frac{1}{2}(y^{-1}-\frac{m_{1}+m_{2}}{2})h_{i}, \
y^{-k}(y^{-1}-m_{1})(y^{-1}-m_{2})e_{\alpha},
$$
$$y^{-k}(y^{-1}-m_{1})(y^{-1}-m_{2})e_{-\alpha},\ \frac{1}{2}(y^{-1}-m_{1})(y^{-1}-m_{2})h_{i},\},\quad
k=1,2,\ldots
$$

\begin{remark}
\label{remark-2}
Complementary Lagrangian subalgebras $W$ are not unique. The corresponding $r$-matrices have the form $r(x,y)=r_{A_{i}}(x,y)+p(x,y)$,
where
$$
p(x,y)\in(\g\otimes \g)[x,y]\quad {\rm and}\quad p(x,y)+p_{21}(y,x)=0.
$$
Later we will present combinatorial data describing all $W$.
\end{remark}
\begin{remark}
For the case $\g=\mathfrak{sl}_2$, the $r$-matrix
$$
r_{A4,m_{1},m_{2}}(x,y)
=\frac{1-(m_{1}+m_{2})u+m_{1}m_{2}uv}{y-x}\Omega-r_{m_{1},m_{2}}
$$
is the classical limit
of the quantum $R$-matrix considered in \cite{FLM} in connection with exactly solvable stochastic processes.
\end{remark}

Now we discuss the two remaining cases
\begin{itemize}
\itemsep-3pt
\item
$D(\g[[x]],\overline{\delta})=\g((x))\oplus \g$ (case B),
\item
$D(\g[[x]],\overline{\delta})=\g((x))\oplus \g[\varepsilon]$, where $\varepsilon^{2}=0$ (case C).
\end{itemize}
We remind the reader that the canonical form in case B is given by the following formula:
$$
(f_{1}(x)+g_{1}|\ f_{2}(x)+g_{2})={t}({K'}(f_{1},f_{2}))-{K}(g_{1},g_{2}),
$$
where the trace $t:\mathbb{C}((x))\longrightarrow \mathbb{C}$ is defined by the formulas:
$t(x^{k})=0\ for\ k\geq 1$, ${t}(1)=1$ and
${t}(x^{-k})=b_{k}$.

Let us denote the series $1+b_{1}x+b_{2}x^{2}+.....$ by $b(x)$. The following results can be proved similarly to that of the case A.

\begin{theorem}
\label{thm-B}
Assume that
$$
\sigma(x)(W(\mathbb{C}[[x^{-1}]]\oplus \mathbb{C}))\subset \mathbb{O}_{\alpha_{i}}\oplus \g.
$$
Then $\frac{1}{b(x)}$ is a polynomial of degree at most 1 if $k_{i}=1$  and $b(x)=1$ if $k_{i}>1$.
\end{theorem}

\begin{corollary}
\label{cor-B}
If $D(\g[[x]],\overline{\delta})=\g((x))\oplus \g$, then $D(\g[x],\delta)=\g[x,x^{-1}]\oplus \g$ as Lie algebras.
The canonical form on $D(\g[x],\delta)$ in this case is given by the formula
$$
(f_{1}(x)+g_{1}|\ f_{2}(x)+g_{2})=Res_{x=0}({K'}(f_{1},f_{2})x^{-1}b(x))-{K}(g_{1},g_{2}).
$$
\end{corollary}
Up to automprphism $\gamma:\mathbb{C}[x,x^{-1}]\longrightarrow \mathbb{C}[x,x^{-1}]$ given by $\gamma(x)=cx$, $c\in \mathbb{C}$,
we have two possibilities $b_{1}(x)=1\ (case\ B1)$ and $b_{2}(x)=\frac{1}{1-x}\ (case\ B2)$.

\begin{remark}\label{remark-B}
At this point we note that we need to present $W_{1}$ and $W_{2}$, which are
Lagrangian subalgebras of $\g[x,x^{-1}]\oplus \g$ transversal to
$\g[x]\subset \g[x,x^{-1}]\oplus \g$ with respect to the canonical forms
determined by $b_{1}(x)$ and $b_{2}(x)$.
\end{remark}
In the case $B1$ we can choose
$$
W_{1}={\rm
Span}\{x^{-1}\g[x^{-1}],(e_{\alpha},0),(0,e_{-\alpha}),(h_{i},-h_{i})\}.
$$
In the case B2, it is easy to verify that
$$
W_{2}={\rm Span}\{(1-x^{-1})\g[x^{-1}],(e_{\alpha},0),(0,e_{-\alpha}),(h_{i},-h_{i})\}$$
satisfies all the conditions we need.
The corresponding $r$-matrices are:
$$
r_{B1}(x,y)=\frac{x}{y-x}\Omega+r_{DJ},\quad
r_{B2}(x,y)=\frac{x(1-y)}{y-x}\Omega+r_{DJ}.
$$
Now let us treat the last case C. We have:
$$
D(\g[[x]],\overline{\delta})=\g((x))\oplus \g[\varepsilon],\quad \varepsilon^{2}=0.
$$

The canonical form in the case C is given by the following formula
\begin{align*}
(f_{1}(x)+g_{1}\varepsilon+h_{1}&|\ f_{2}(x)+g_{2}\varepsilon+h_{2})\\
&=t({K'}(f_{1},f_{2})-{K}(g_{1},h_{2})-{K}(h_{1},g_{2})
-c_{1}{K}(h_{1},h_{2})),
\end{align*}
where $t:\mathbb{C}((x))\longrightarrow \mathbb{C}$
is defined in the following way :
$t(x^{k})=0$ for $k\geq 2$, $t(x)=1$, $t(1)=c_{1}$, $t(x^{-k})=c_{k+1}$.

As before, let us denote series $1+c_{1}x+c_{2}x^{2}+\ldots$ by $c(x)$. Then we have the following results.
\begin{theorem}
\label{theorem-C}
Assume that
$$
\sigma(x)(W(\mathbb{C}[[x^{-1}]]\oplus \mathbb{C}))\subset \mathbb{O}_{\alpha_{i}}\oplus \g[\varepsilon].
$$
Then
\begin{itemize}
\itemsep-1pt
\item
$c(x)=1$ if $k_i=1$.
\item
$\mathbb{O}_{\alpha_{i}}\oplus \g[\varepsilon]$ does not contain Lagrangian subalgebras if $k_i>1$.
\end{itemize}
\end{theorem}
This theorem implies that the case C coincides with the so-called $4^{th}$
structure considered in \cite{SY}.
In particular, $W=\g[x^{-1}]+\varepsilon \g$ is a Lagrangian subalgebra transversal to $\g[x]$ embedded into $\g[x,x^{-1}]\oplus \g[\varepsilon]$
as follows:
$$
a_0+a_1 x+...+a_n x^n\longrightarrow (a_{0}+a_{1}x+...+a_{n}x^{n})\oplus(a_{0}+a_{1}\varepsilon).
$$
The corresponding $r$-matrix is
$$
r_{C}(x,y)=\frac{xy}{y-x}\Omega.
$$

\section{Description of all Lie bialgebra structures on current polynomial Lie algebras}

The aim of this section is to describe all Lagrangian subalgebras $W$ of $D(\g[x],\delta)$ such that $D(\g[x],\delta)=\g[x]\oplus W$,
direct sum of vector spaces.

It will be done along with classification of solutions of the classical Yang--Baxter equation of certain types.
Throughout this section we assume that $W\oplus \g[x]=D(\g[x],\delta)$.

First of all, we have proved that $W\subset x^{N}\g[x^{-1}]$.
Furthermore, we know that applying some $\sigma(x)\in Aut_{\mathbb{C}[x]}(\g[x])$,
we can achieve that
\begin{itemize}
\itemsep-3pt \item $\sigma(x)W\subset \mathbb{O}_{\alpha_{i}}$ in the case
A, \item $\sigma(x)W\subset \mathbb{O}_{\alpha_{i}}\oplus \g$ in the case B,
\item $\sigma(x)W\subset \mathbb{O}_{\alpha_{i}}\oplus \g[\varepsilon]$ in
the case C.
\end{itemize}

Let $\mathfrak{p}^{-}_{\alpha_{i}}\subset g$ be the parabolic subalgebra
corresponding to $\alpha_{i}$, i.e. is generated by all $e_{-\alpha}$ and
those $e_{\alpha}$, which do not contain $\alpha_{i}$ in their simple root
decomposition.
\begin{definition}\label{definition-C}
The data $F(\alpha_{i},k_{i},L,B)$ consists of a subalgebra $L\subset \g$
such that $L+\mathfrak{p}^{-}_{\alpha_{i}}=\g$ and a $2$-cocycle $B$ on $L$
such that $B$ is nondegenerate on $L\cap \mathfrak{p}^{-}_{\alpha_{i}}$.
\end{definition}

It was proved in \cite{S1}, \cite{S2} that in the case A1 there exists a one-to-one  correspondence between the sets
$\{W\subset\mathbb{O}_{\alpha_{i}},k_{i}=1\}$ and $\{F(\alpha_{i},1,L,B)\}$. In the sequel we will use the following notation
$\{W\subset\mathbb{O}_{\alpha_{i}},k_{i}=1\}\stackrel{1-1}{\longleftrightarrow}\{F(\alpha_{i},1,L,B)\}$.

\begin{remark}
\label{remark-D}
In the case A1 the set $\{W\subset\mathbb{O}_{\alpha_{i}},k_{i}=3\}$ has a similar but more complicated description and not much is known for other $k_{i}$.
\end{remark}

It turns out that in the cases A3 and C the corresponding Lagrangian subalgebras transversal to $\g[x]$
can be completely described by the data $\{F(\alpha_{i},1,L,B)\} $.

\begin{theorem}
\label{theorem-D}
1. Case A3:

(a) $\{W\subset\mathbb{O}_{\alpha_{i}},k_{i}=1\}\stackrel{1-1}{\longleftrightarrow}\{F(\alpha_{i},1,L,B)\}, $

(b) the set $\{W\subset\mathbb{O}_{\alpha_{i}},k_{i}\geq 2\}$ is empty.

2. Case C:

(a) $\{W\subset\mathbb{O}_{\alpha_{i}}\oplus \g[\varepsilon], k_{i}=1\}
\stackrel{1-1}{\longleftrightarrow}\{F(\alpha_{i},1,L,B)\},$

(b) the set $\{W\subset\mathbb{O}_{\alpha_{i}}\oplus \g[\varepsilon], k_{i}\geq 2\}$ is empty.
\end{theorem}

\begin{proof}
It will be explained later why Cases A3 and C have one and the same description. Then Case C(a) was considered in
\cite{PY},  Case A3(b) follows from Theorem \ref{theorem-ma}.
\end{proof}

For a description of the Lagrangian subalgebras transversal to $\g[x]$ in the remaining cases
we will define new data $BD(\alpha_{i},k_{i},\Gamma_{1},\Gamma_{2},\tau,s)$.

\begin{definition}\label{definition-E}
$BD(\alpha_{i},k_{i},\Gamma_{1},\Gamma_{2},\tau,s)$ consists of:
\begin{enumerate}
\itemsep-3pt \item $\alpha_{i}\in \widehat{\Gamma}$, which is the extended
Dynkin diagram of $\g$, \item $\Gamma_{1},\Gamma_{2}\subset
\widehat{\Gamma}$ which are such that $\alpha_{i}\notin\Gamma_{1}$ and
$\alpha_{0}\notin\Gamma_{2}$. \item $\tau:\Gamma_{1}\longrightarrow
\Gamma_{2}$ is admissible in sense of Belavin-Drinfeld, i.e.
$\langle\tau(\alpha),\tau(\beta)\rangle=\langle\alpha,\beta\rangle$ for all
$\alpha,\beta\in \Gamma_{1}$ and $\tau^{n}$ is not defined for large $n$.
\item $s\in \Lambda^{2}V$, where $V\subset \h$ is defined by the linear
system $(\alpha-\tau(\alpha))(h)=0$, $\alpha\in\Gamma_1$.
\end{enumerate}
\end{definition}

The following theorem can be deduced from \cite{PS}

\begin{theorem}
\label{theorem-E}
In the case B1, $\{W\subset \mathbb{O}_{\alpha_{i}}\oplus \g\}\stackrel{1-1}{\longleftrightarrow}\{BD(\alpha_{i},k_{i},\Gamma_{1},\Gamma_{2},\tau,s)\}$.
\end{theorem}
Later we will explain why this result implies

\begin{theorem}
\label{theorem-F}
The set $\{W\subset \mathbb{O}_{\alpha_{i}}\}\stackrel{1-1}{\longleftrightarrow}
\{BD(\alpha_{i},k_{i},\Gamma_{1},\Gamma_{2},\tau,s)\}$ in the case A2.
\end{theorem}

In the remaining cases A4 and B2, the results are analogous to the theorems
above.

\begin{theorem}
\label{theorem-G}
Both sets $\{W\subset\mathbb{O}_{\alpha_{i}},k_{i}=1\}$ (case A4) and $\{W\subset\mathbb{O}_{\alpha_{i}}\oplus \g,k_{i}=1\}$ (case B2)
are in a one-to-one correspondence with the set $\{BD(\alpha_{i},1,\Gamma_{1},\Gamma_{2},\tau,s)\}$.
The sets $\{W\subset\mathbb{O}_{\alpha_{i}}\}$ (case A4) and $\{W\subset\mathbb{O}_{\alpha_{i}}\oplus \g\}$ (Case B2) are empty if $k_{i}>1$.
\end{theorem}

\section{Connection with solutions of the classical Yang Baxter equation}

We already know that the classical $r$-matrices $r_{X_{i}} (x,y)$ define Lie bialgebras on $\g[x]$.
Clearly, if $r(x,y)=r_{X_{i}}(x,y)+p(x,y)$ satisfies the classical Yang-Baxter equation,
then $r(x,y)$ defines a Lie bialgebra on $\g[x]$. Here $X=A,B,C$ and $p(x,y)\in (\g\otimes \g)[x,y]$ is such that $p(x,y)=-p^{21}(y,x)$.

It is well known that the classical doubles constructed from $r_{X_{i}} (x,y)$ and $r(x,y)=r_{X_{i}}(x,y)+p(x,y)$
are isomorphic as Lie algebras with canonical forms (for proofs see \cite {S1} and \cite{KPSST}).

Therefore, the following result is true:

\begin{theorem}
\label{theorem-I} Solutions of the classical Yang-Baxter equation of the
form
$$
r(x,y)=r_{X_{i}}(x,y)+p(x,y)
$$
are in a 1-1 correspondence with the Lagrangian subalgebras $W$ of the respective double satisfying conditions
\begin{enumerate}
\itemsep-3pt
\item
$W\oplus \g[x]=D(\g[x],\delta)$ (cases A, B, C).
\item
$W\subset x^{N}\g[x^{-1}]$ (case A),
$W\subset x^{N}\g[x^{-1}]\oplus \g$ (case B) and
$W\subset x^{N}\g[x^{-1}]\oplus \g[\varepsilon]$ (case C).
\end{enumerate}
\end{theorem}

\begin{corollary}
\label{corollary-I}
In the following pairs of cases,
there is a one-to-one  correspondence between the corresponding sets
of Lagrangian subalgebras of $D(\g[x])$ transversal to $\g[x]$:
$A2$ and $B1$,
$A3$ and $C$,
$A4$ and $B2$.
\end{corollary}

\begin{proof}
We start the proof with the following observation:
In each case there exists  automorphisms of $\g[x]$ of the form
$\tau_i: x\to px+q$, $p,q,s_i\in\mathbb{C}$
such that
$$
(\tau_1 \otimes\tau_1 )(r_{B1})=s_1 r_{A2} ,\quad
(\tau_2 \otimes\tau_2 )(r_{C})=s_2 r_{A3} ,\quad
(\tau_1 \otimes\tau_1 )(r_{B2})=s_3 r_{A4}
$$
(it immediately follows from Corollary \ref{corollary-3} and Remark \ref{remark-B}). Therefore, the correspondence between the three cases exists.

Let us explain how the corresponding Lagrangian subalgebras relate to each other.
Consider for instance the cases A3 and C. Assume that  $\alpha_i$ has
coefficient $1$ in  decomposition $\alpha_{max}=\sum k_n \alpha_n$. We can also assume
that the Lagrangian subalgebra $W_C$ is contained in $\mathbb{O}_{\alpha_i} \oplus \g[\varepsilon]$.
Then $W_C$ contains the orthogonal complement to $\mathbb{O}_{\alpha_i} \oplus \g[\varepsilon]$.
The latter is $\mathbb{O}_{\alpha_i} $. Therefore, $W_C$ is uniquely defined by its image $\widehat{W_C}$
under the canonical projection
$$
\frac{\mathbb{O}_{\alpha_i} \oplus \g[\varepsilon]}{\mathbb{O}_{\alpha_i} }\to \g[\varepsilon]
$$
Since $W_C$ is transversal to $\g[x]$, it is easy to see that
$\widehat{W_C}$ should be transversal to the image of $\g[x]\cap
({\mathbb{O}_{\alpha_i} \oplus \g[\varepsilon]})$ under the same projection.
This image is isomorphic to $\mathfrak{p}^{-}_{\alpha_{i}}
+\varepsilon(\mathfrak{p}^{-}_{\alpha_{i}})^{\perp}$, where
$(\mathfrak{p}^{-}_{\alpha_{i}})^{\perp}$ is the orthogonal complement to
$\mathfrak{p}^{-}_{\alpha_{i}}$ in $\g$ with respect to the Killing form.

Let us turn to the case A3. We are looking for a Lagrangian subalgebra $W_{A3}\subset \mathbb{O}_{\alpha_i}\cap \g[x,x^{-1}]$
transversal to $\g[x]$. Then $W_{A3}$ contains the orthogonal complement to $\mathbb{O}_{\alpha_i}\cap \g[x,x^{-1}]$, which
is $x^{-2}(1-x)^2 (\mathbb{O}_{\alpha_i}\cap \g[x,x^{-1}])$ because $\alpha_i $ has coefficient $1$ in the decomposition of
the maximal root, and it was proved in ~\cite{S2} that
$$
\mathbb{O}_{\alpha_i}\cap \g[x,x^{-1}]={\rm Ad} (X(x)) (\mathbb{O}_{\alpha_0}\cap \g[x,x^{-1}] ),
$$
where $X(x)\in {\rm Ad} (\g[x,x^{-1}])$. It follows that
$$
\frac{\mathbb{O}_{\alpha_i}\cap \g[x,x^{-1}]}{x^{-2}(1-x^2)(\mathbb{O}_{\alpha_i}\cap \g[x,x^{-1}])}\cong \g[\gamma],
$$
where $\gamma=x^{-1} -1$ and $\gamma^2 =0$. Again we see that $W_{A3}$ is
uniquely defined by its image $\widehat{W_{A3}}$ in $\g[\gamma]$ and
$\widehat{W_{A3}}$ should be transversal to the image of $\g[x]\cap
{\mathbb{O}_{\alpha_i} }$ in $\g[\gamma]$. The latter is again
$\mathfrak{p}^{-}_{\alpha_{i}} +\gamma
(\mathfrak{p}^{-}_{\alpha_{i}})^{\perp}$, what completes the proof of the
theorem for the cases A3 and C. Two other cases can be considered in a
similar way.
\end{proof}

\section{Remarks about quantization}

It was proved in \cite{H} that if $(L_{1},\delta_{1})$ and $(L_{2},\delta_{2})$ are such that there exists an isometry
$\varphi:D(L_{1},\delta_{1})\longrightarrow {D}(L_{1},\delta_{2})$ (with respect to the canonical forms)
satisfying $\varphi(L_{1})=L_{2}$, then we can find quantizations of $(L_{1},\delta_{1})$
and $(L_{2},\delta_{2})$ denoted by $(H_{1},\Delta_{1})$ and $(H_{2},\Delta_{2})$
such that $H_{1}$ and $H_{2}$ are isomorphic as algebras. Moreover, if
$\overline{\varphi}:H_{1}\longrightarrow H_{2}$ is this isomorphism,
then
$$
(\overline{\varphi}^{-1}\otimes \overline{\varphi}^{-1})\Delta_{2}(\overline{\varphi}(a))=F\Delta_{1}(a)F^{-1},
$$
where $F$ satisfies the so called cocycle equation
$$
F_{12}(\Delta_{1}\otimes id)(F)=F_{23}(id\otimes \Delta_{1})(F).
$$
We call $F$ a quantum twist.

Hence, we see that all Lie bialgebras related to one and the same case
$X_{i}$ ($X=A,B,C$) can be quantized in such a way that they will be
isomorphic as algebras and their co-algebra structures will differ by a
quantum twist $F$.
\begin{itemize}
\itemsep-3pt
\item
Quantization of $A1$ is called Yangian \cite{DR}.
\item
Quantization of $B1$ is $U_{q}(\g[x])$ \cite{PS}.
\item
One can show that quantization of $C1$ is $U(\g)\otimes RDY(\g)$, where $RDY(\g)$ is the restricted dual Hopf algebra to the Yangian $Y(\g)$.
\item
Quantization of $A2$ is the so-called Drinfeldian, $Drin(\g)$ (see \cite{T}).
\item
Quantizations of the cases $A3$, $B2$, $A_{4,m_{1},m_{2}}$ are not known yet.
\end{itemize}


\small

\noindent F.M.: Department of Mathematics, University of Zaragoza\\
50009 Zaragoza, Spain\\
e-mail: {\small \tt fmontane@unizar.es}

\medskip

\noindent A.S.: Department of Mathematics, University of Gothenburg,\\
SE-412 96 G\"oteborg, Sweden\\
e-mail: {\small \tt alexander.stolin@gu.se}

\medskip

\noindent E.Z.:  Department of Mathematics,UCSD\\
9500 Gilman Drive\\
La Jolla, CA 92093-0112, USA\\
e-mail: {\small \tt ezelmano@math.ucsd.edu}


\begin{thebibliography}{10}

\bibitem{BD} Belavin, A., Drinfeld, V. \textit{Triangle equations and simple
algebras}. Math. Phys. Rev. 4(1984), 96-165. \bibitem{BZ} Benkart, G.,
Zelmanov, E. \textit{Lie algebras graded by the finite root system and
intersection matrix algebras}. Invent. Math. 126(1996), 1-145. \bibitem{B}
Benkart, G. On inner ideals and ad-nilpotent elements of Lie algebras .
Trans. Amer. Math. Soc. 232 (1977) 61-81. \bibitem{De} Delorme, P.
\textit{Classification des triples de manin pour les algebras de Lie
reductivs complexes}, T. Algebra 246(2001), 97-174. \bibitem{DR} Drinfeld,
V. \textit{Quantum Groups}. Proceedings of the International Congress of
Mathematicians (Berkeley, 1986), 798-820. \bibitem{EK} Etigof, P., Kazhdan,
D. \textit{Quantization of Lie Bialgebras. I--V.} Selecta Math. (N. S.)
2(1996), no. 1, 1--41, 4(1998), no. 2, 213--231, 233--269, 6(2000), no. 1,
79--104, 105--130. \bibitem{EK1} Etigof, P., Kazhdan, D.
\textit{Quantization of Lie Bialgebras VI. Quantization of generalized
Kac--Moody algebras.} Transform. Groups 13(2008), no. 3--4, 527--539.
\bibitem{FLM} Furtlechner, C., Lesgouttes, J.-M., Samsonov, M.
\textit{Exactly Solvable Stochastic Processes for Traffic Modelling.}
Preprint of the French National Institute for Research in Computer Science
and Control (INRIA). 2010. \bibitem{H} Halbout, G. \textit{Formality theorem
for Lie Bialgebras and quantization of twists and coboundary r-matrices.}
Adv. Math. 207(2007), 617-633. \bibitem{J} Jacobson, Nathan, Lie algebras,
Republication of the 1962 original. Dover Publications, Inc., New York,
1979. \bibitem{KPSST} Khoroshkin, S., Pop, I., Samsonov, M., Stolin, A.,
Tolstoy, V. \textit{On some Lie Bialgebra structures on polynomial algebras
and their quantization.} Comm. Math. Phys. 282(2008), 625-662. \bibitem{PS}
Pop, I., Stolin, A. \textit{Lagrangian subalgebras and quasi-trigonometric
r-matrices}. Lett. Math. Phys. 85 (2008), 249-262. \bibitem{PY} Pop, I.,
Yermolova-Magnusson, J. \textit{New r-matrices for Lie Bialgebra Structures
over Polynomials}. ArXiv: 0910.4286 \bibitem{S1} Stolin, A. \textit{On
rational solutions of Yang-Baxter equation for sl(n).} Math. Scand. 69
(1991), 57-80. \bibitem{S2} Stolin, A. \textit{Rational solutions of CYBE.
Maximal orders in loop algebras}. Comm. Math. Phys. 143(1991). \bibitem{S3}
Stolin, A. \textit{A geometrical approach to rational solutions of the
classical Yang-Baxter equation.} Proceedings of the 2nd Gauss Symposium.
Conference A: Mathematics and Theoretical Physics (Munich, 1993), \textit{
347--357, Sympos. Gaussiana,} de Gruyter, Berlin, 1995. \bibitem{SY} Stolin,
A., Yermolova-Magnusson, J.\textit{ The 4th structure}. (Czech. J. Phys.)
56(2006), 1296-1297.

\bibitem{T} Tolstoy, V. \textit{From quantum affine Kac-Moody algebras to Drinfeldians and Yangians}.  Contemp. Math. 343, 349-370.

\end{thebibliography}
\end{document}